\documentclass[a4paper,reqno]{amsart}

\usepackage[all]{xy}
\usepackage{color}
\usepackage{mathrsfs}
\usepackage{hyperref}
\usepackage{enumerate}
\usepackage{amsmath}
\usepackage{amsthm}
\usepackage{amssymb}
\usepackage{amscd}
\usepackage{graphicx}
\usepackage{epsfig}
\usepackage[english]{babel}
\usepackage[stable]{footmisc}

\newtheorem{proposition}{Proposition}
\newtheorem{lemma}{Lemma}
\newtheorem{theorem}{Theorem}
\newtheorem{corollary}{Corollary}
\newtheorem{definition}{Definition}
\newtheorem{example}{Example}
\newtheorem{remark}{Remark}

\begin{document}

\author{David Mart\'inez Torres}

\email{dfmtorres@gmail.com}

\address{Centro de An\'{a}lise Matem\'{a}tica, Geometria e Sistemas
Din\^{a}micos,
Departamento de Matem\'atica, Instituto Superior T\'ecnico, Av. Rovisco Pais,
1049-001
Lisboa, Portugal}

\title[2-calibrated foliations]{Codimension one foliations calibrated by non-degenerate closed 2-forms}

\begin{abstract}
 A  class of codimension one foliations has been  recently introduced by
 imposing a natural compatibility condition with a closed maximally non-degenerate 2-form. 
  In this paper we study for such foliations the information captured by a Donaldson type submanifold. In particular we deduce 
that their leaf spaces are homeomorphic to leaf spaces of 3-dimensional taut foliations. We also introduce surgery
constructions to show that this class of foliations is broad enough. Our techniques come mainly from symplectic geometry.
\end{abstract}
\maketitle

%
%

\section{Introduction and Statement of Main Results}

Codimension one foliations are too large a class of structures to
obtain strong structure theorems for them. According to a theorem of
Thurston \cite{Th76} a closed manifold admits a codimension one
foliation if and only if its Euler characteristic  vanishes. In
order to draw significant results it is necessary to assume the
existence of other structures compatible with the foliation. 

From the point of view of symplectic geometry it is natural to consider the following class of codimension
 one foliations:
\begin{definition} \cite{IM04a}\label{def:main} A codimension one foliation
$\mathcal{F}$ of $M^{2n+1}$ is
said to be 2-calibrated if there exists a closed $2$-form $\omega$ such
that  ${\omega_{\mathcal{F}}}^{n}$ is no-where vanishing (we also say that $\omega^n$ is no-where vanishing on $\mathcal{F}$).

The 2-calibrated foliation is said to be integral  if $[\omega]\in H^2(M;\mathbb{Z})$.
\end{definition}
The notation ${\omega_{\mathcal{F}}}^{n}$ in definition \ref{def:main} stands for the restriction of $\omega^n$ 
to the leaves of $\mathcal{F}$. 
We will be using the subscript $\mathcal{F}$ 
(respectively $W$, if $W$ is a
submanifold of $M$) to denote the restriction of a form,
connection, etc,  to the leaves of $\mathcal{F}$ (respectively to $W$).  In what follows the  manifolds will
always be closed and oriented, the codimension one foliations
co-oriented and all the structures and maps smooth.

In the next paragraphs we are going to describe how the 2-calibrated condition appears naturally 
when looking at the problem of
constructing submanifolds transverse to a codimension one foliation.

Recall that a codimension one foliation  $\mathcal{F}$ is said to be taut if every leaf meets a  transverse 1-cycle. 
Tautness in codimension one can be characterized in several ways 
using forms, metrics and currents \cite{Su76,Ru79,HL82}. The characterization we are interested in, says
 that a rank $p$ codimension one foliation $\mathcal{F}$ is taut if and only if there exists a closed $p$-form  $\xi$ 
 no-where vanishing on $\mathcal{F}$
(and furthermore according to proposition 2.7 in \cite{HL82}, it is possible to construct a metric
 $g$ so that $\xi$ is a calibration for $(M,\mathcal{F})$). Note in particular that a 2-calibrated 
foliation $(M,\mathcal{F},\omega)$ is always taut,
 since $\xi:=\omega^n$ is no-where vanishing on $\mathcal{F}$. In dimension three 2-calibrated foliations are the
same as taut foliations.

Let us analyze one direction of the aforementioned characterization: the existence of 
a closed $p$-form whose restriction to each leaf is a volume form,  is
 equivalent to a reduction of the structural pseudogroup of
$(M,\mathcal{F})$ to
$\mathrm{Vol}(\mathbb{R}^{p},\Xi_{\mathbb{R}^{p}})\times \mathrm{Diff}(\mathbb{R})$, where
\[\Xi_{\mathbb{R}^{p}}:=dx_1\wedge \cdots\wedge dx_p,\]
$x_1,\dots,x_p$ are coordinates on $\mathbb{R}^p$, and $\mathrm{Vol}(\mathbb{R}^{p},\Xi_{\mathbb{R}^{p}})$
 (respectively  $\mathrm{Diff}(\mathbb{R})$)
is the pseudogroup of
local diffeomorphisms of $\mathbb{R}^{p}$ (respectively $\mathbb{R}$) preserving the
volume form $\Xi_{\mathbb{R}^{p}}$.
 Let $U$ be any open  subset of a leaf of $\mathcal{F}$. Poincar\'e recurrence theorem  implies that 
the flow of any vector field spanning $\mathrm{ker}\xi$ defines a
first return map from $U'\subset U$ to $U''\subset U$. A straightforward consequence is that
closed transverse 1-cycles through any given $x\in M$ can be constructed by slightly
deflecting integral curves of $\mathrm{ker}\xi$.

 The first return map
belongs to the pseudogroup $\mathrm{Vol}(\mathbb{R}^{p},\Xi_{\mathbb{R}^{p}})$. If $p=2$, that is, if we
have a taut foliation on a 3-manifold, then under certain
circumstances we can deduce interesting geometric information about
the existence of more closed orbits (Poincar\'e-Birkhoff theorem).
If $p>2$ we have little geometric control on the return map because
assuming for simplicity that $U'$ and $U''$ are diffeomorphic to a
ball, the only invariant is the total volume (\cite{GS80}, theorem 1).  Therefore problems such as the
existence of transverse submanifolds of dimension bigger than one seem
difficult to attack.

It has been known for some time that the right setting to obtain higher dimensional
generalizations of Poincar\'e-Birkhoff theorem is not volume geometry but symplectic geometry
 (\cite{HZ94}, chapter 6; \cite{MS98}, chapter IV).
It can be checked (see section \ref{sec:def}) that the existence of a closed 2-form $\omega$ which makes the leaves of   
 $(M,\mathcal{F})$ symplectic manifolds, amounts to a 
 reduction of the structural pseudogroup of  $(M,\mathcal{F})$ to
$\text{Symp}(\mathbb{R}^{2n},\Omega_{\mathbb{R}^{2n}})\times \text{Diff}(\mathbb{R})$,
where $\text{Symp}(\mathbb{R}^{2n},\Omega_{\mathbb{R}^{2n}})$ is the pseudogroup of local diffeomorphisms 
of $\mathbb{R}^{2n}$ preserving
the  standard symplectic form
\[
\Omega_{{\mathbb{R}^{2n}}}:=\sum_{i=1}^ndx_i\wedge dy_i.
\]
 Thus, return maps associated to the flow of vector fields generating
$\mathrm{ker}\omega$ belong to
$\mathrm{Symp}(\mathbb{R}^{2n},\Omega_{\mathbb{R}^{2n}})$.
Symplectomorphisms are much more rigid than 
transformations preserving the volume form
$\Omega_{\mathbb{R}^{2n}}^n=n!\Xi_{\mathbb{R}^{2n}}$.
They preserve the symplectic invariants of  subsets
of $\mathbb{R}^{2n}$, so for example these cannot be squeezed along
symplectic 2-planes  (\cite{HZ94}, chapters 2 and 3; \cite{MS98}, section 12). Naively, one might try to construct transverse
3-manifolds  by choosing tiny 2-dimensional
symplectic pieces $\Sigma$ inside a leaf, whose image by  the
first return map
is a small 2-dimensional symplectic manifold
that can be isotoped  to
$\Sigma$ through symplectic surfaces. The isotopy would be
used to connect both symplectic surfaces in nearby leaves, and thus get a piece
of transverse 3-dimensional taut foliation. Of
course this idea seems difficult to be carried out because 
different pieces should be combined to construct a closed 3-manifold. However,
 it provides some insight on why  2-calibrated foliations are expected to have embedded 3-dimensional taut
foliations. 

In \cite{IM04a}, corollary 1.2, it was proved that for  any  2-calibrated foliation  $(M,\mathcal{F},\omega)$
there exists an embedding  of a 3-dimensional submanifold $W^3\hookrightarrow M$, such that
$W^3$ is transverse to $\mathcal{F}$ and  $\omega_W$ is no-where vanishing on $\mathcal{F}_W$;
 the 3-dimensional submanifold $W^3$, 
which inherits a taut
foliation, is a 
Donaldson type submanifold \cite{Do96,Au98}. Its existence is
 an elementary consequence of the extension to 2-calibrated foliations 
of the approximately holomorphic techniques for symplectic manifolds
introduced by Donaldson \cite{Do96}.

\subsection{Statement of results}

Let  $(M,\mathcal{F},\omega)$ be a 2-calibrated foliation and let $W\hookrightarrow (M,\mathcal{F})$
 be a 3-dimensional Donaldson type submanifold. 
 In this paper we are mainly concerned with finding out
which  properties of $(M,\mathcal{F})$ are captured
 by $W$. 

If $F$ is a compact leaf of $(M,\mathcal{F},\omega)$,
an appropriate version of  the Lefschetz hyperplane theorem  (\cite{Do96}, proposition 39) asserts that $W\cap F$
is connected. A codimension one foliation  $(M,\mathcal{F})$ has non-compact leaves
unless it is a fibration over the circle (a  mapping torus). If $F$ is a non-compact leaf then
 describing global properties of
$W\cap F$ seems very difficult.  Our main result is a rather
surprising and counterintuitive global property of such intersections for
appropriate Donaldson type  submanifolds.

\begin{theorem}\label{thm:lefhyp} Let $(M,\mathcal{F},\omega)$ be a  2-calibrated foliation.
Then there exists Donaldson type  submanifolds  $W^3\hookrightarrow
(M,\mathcal{F})$, such that for every leaf $F$ of $\mathcal{F}$ the intersection $W\cap F$ is
connected.
\end{theorem}

\begin{remark} Any integral 2-calibrated foliation $(M,\mathcal{F},\omega)$ admits embeddings in complex projective
spaces $\mathbb{C}\mathbb{P}^N$ of large dimension,  with the property that  the
 ambient Fubini-Study symplectic form restricts to a multiple of $\omega$ (\cite{IM04a}, corollary 1.3).
 The 3-dimensional transverse submanifolds in theorem \ref{thm:lefhyp}
can be arranged to appear as intersections of $M\subset \mathbb{C}\mathbb{P}^N$ with 
appropriate projective subspaces. Theorem \ref{thm:lefhyp}
 should be understood as a leafwise Lefschetz hyperplane type
result for  $\pi_0$.
\end{remark}

An important consequence of theorem \ref{thm:lefhyp} is the following result:
\begin{theorem}\label{thm:leafhomeo} Let $(M,\mathcal{F},\omega)$ be a  2-calibrated foliation.
Then there exists a 3-dimensional embedded taut foliation such that the inclusion
$(W^3,\mathcal{F}_W)\hookrightarrow (M,\mathcal{F})$ descends to a homeomorphism of leaf
spaces $W/\mathcal{F}_W\rightarrow M/\mathcal{F}$.

Thus, leaf spaces of 2-calibrated foliations are no more complicated than those
of 3-dimensional taut foliations.
\end{theorem}

A second goal of this paper is showing that 2-calibrated foliations are a
broad enough class of foliations. In this respect 
there are three basic families of 2-calibrated foliations: products,  cosymplectic
foliations and symplectic bundle foliations.

In a product we cross a 2-calibrated foliation
-typically a 3-dimensional taut foliation- with a (non-trivial) symplectic manifold, and put
the product foliation and the obvious closed 2-form.

A cosymplectic foliation is a triple $(M,\alpha,\omega)$, where $\alpha$ is a no-where vanishing closed 1-form and 
$(M,\mathrm{ker}\alpha,\omega)$ is a 2-calibrated foliation. 

A bundle foliation with fiber $S^1$ is by definition an $S^1$-fiber bundle $\pi\colon M\rightarrow X$ endowed
 with a codimension one foliation $\mathcal{F}$ transverse to the
fibers. If the base space admits a symplectic form $\sigma$, then $(M,\mathcal{F},\pi^*\sigma)$ is a 2-calibrated foliation which
 we refer to as a symplectic bundle foliation.

The second topic of this paper concerns the introduction of two surgery
constructions for  2-calibrated foliations: normal connected sum  and generalized Dehn surgery
or Lagrangian surgery.  Using surgery we have obtained the following result:
\begin{proposition}\label{pro:newex}
There exist 2-calibrated foliations 
(of dimension bigger than three) which are neither products, nor cosymplectic foliations, nor symplectic bundle foliations. 
\end{proposition}

The paper is organized as follows. In section \ref{sec:def} we
introduce definitions and basic facts on 2-calibrated foliations, and address
their relation to regular Poisson
structures.

Section \ref{sec:nsum} describes how to adapt the normal connected sum for
symplectic and Poisson manifolds to integral 2-calibrated foliations; this is the surgery
used to prove proposition \ref{pro:newex}.

In section \ref{sec:dsurg} we present a surgery  based on generalized Dehn twists. Generalized Dehn surgery is the
natural extension to 2-calibrated foliations of positive Dehn surgery along a
curve in a leaf of a 3-dimensional taut foliation $(M^3,\mathcal{F})$. 

It is a classical result of
Lickorish \cite{Li65}
 that positive Dehn surgery along a
curve  $\gamma$ has an alternative description: $\gamma$ carries a canonical framing and therefore it determines
an elementary cobordism from $M^3$ to $M'$, which amounts to attaching a 2-handle to the trivial cobordism $M\times [0,1]$. 
The ``new'' boundary component $M'$ is endowed with a canonical foliation which coincides with positive Dehn
surgery on $(M,\mathcal{F})$ along $\gamma$. 

If  $(M^{2n+1},\mathcal{F},\omega)$ 
is a 2-calibrated foliation, a parametrized Lagrangian $n$-sphere inside a leaf of $\mathcal{F}$ canonically
 determines the attaching of a $(n+1)$-handle. We show that the corresponding elementary
$(2n+2)$-dimensional cobordism admits a symplectic structure, which induces a 2-calibrated foliation
on the new boundary component of the cobordism. We call this
construction Lagrangian surgery. In theorem \ref{thm:equiv} we extend Lickorish' result by proving
that generalized Dehn surgery and Lagrangian surgery produce equivalent
2-calibrated foliations. The importance of this result stems from the fact that
the aforementioned symplectic elementary cobordisms do appear in a natural way
associated to Lefschetz pencil structures.
 As a byproduct we get an application to contact geometry that we have included in an appendix:  
it is a proof of a result
announced by Giroux and Mohsen \cite{GM02}, relating generalized Dehn
surgery along a parametrized Lagrangian sphere $L$ in an open book decomposition
compatible with a contact structure, and Legendrian surgery along
$L$. Results in this section require a fine
analysis of the symplectic monodromy about the singular fiber of the complex quadratic form.

In section \ref{sec:pencils} we prove theorems \ref{thm:lefhyp} and \ref{thm:leafhomeo}.
The main tool are  Lefschetz
pencil structures for $(M,\mathcal{F},\omega)$, which are appropriate analogs of leafwise complex Morse
functions and whose existence is an application of approximately
holomorphic geometry for 2-calibrated foliations.  A regular fiber of a Lefschetz pencil structure is a Donaldson type 
submanifold.
A Lefschetz pencil structure admits a leafwise symplectic connection. Its associated
leafwise symplectic parallel transport is the key ingredient to prove our main theorem
relating the leaf space of any regular fiber of the pencil to the leaf space
of $(M,\mathcal{F},\omega)$. Symplectic parallel transport also allows us to compare the 2-calibrated foliations
induced on
different regular fibers. Namely, in theorem \ref{thm:handle} we show that any two regular fibers
of a Lefschetz pencil structure for $(M,\mathcal{F},\omega)$ are related by a sequence of symplectic
handle attachings along Lagrangian spheres. By the symplectic analog 
of Lickorish's result proved in section \ref{sec:dsurg},  we conclude that  any two regular fibers
of a Lefschetz pencil structure are related by a sequence of generalized Dehn surgeries. We finish the
 section by discussing some open problems.

The author is very grateful to the referee for his/her corrections and numerous suggestions. 

\section{Definitions and basic results}\label{sec:def}
In this section we introduce some basic definitions, results and examples. We also address the relation of 2-calibrated
foliations to Poisson structures.

\begin{definition}\label{def:submnfd} Let $(M,\mathcal{F},\omega)$ be a 2-calibrated
foliation
and let $l\colon N\hookrightarrow M$ be a submanifold. We say that $N$ is a
2-calibrated submanifold if $(N,l^*\mathcal{F},l^*\omega)$ is a 2-calibrated
foliation.
\end{definition}

The definition of a 2-calibrated foliation can be given locally.

\begin{definition}\label{def:localdef} A 2-calibration  for $(M,\mathcal{F})$ is a
reduction
of its structural pseudogroup to
$\text{Symp}(\mathbb{R}^{2n},\Omega_{\mathbb{R}^{2n}})\times \text{Diff}(\mathbb{R})$.
\end{definition}

Definitions \ref{def:main} and \ref{def:localdef}  are
equivalent. A standard Darboux type result (see for example \cite{MS98}, chapter 3, 
for basic material on symplectic geometry) implies that about any point in $M$ there exists a foliated chart with coordinates
$x_1,y_1,\dots,x_n,y_n,t$  (the image of $\mathcal{F}$ in $\mathbb{R}^{2n+1}$ is the foliation
by affine hyperplanes with constant coordinate $t$), such that $\omega$ is the pullback of

\[\omega_{\mathbb{R}^{2n+1}}:=\sum_{i=1}^n dx_i\wedge dy_i.\]

It is clear that on a given manifold 2-calibrated foliations are an open subset of the set of codimension one foliations in
 the $C^0$-topology. More precisely, in the product space  of codimension one foliations and  closed 2-forms,
pairs corresponding to 2-calibrated foliations are an open set in the  $C^0$-topology.

The first examples of 2-calibrated manifolds are 3-dimensional taut foliations.
In this paper we are concerned with higher dimensional
2-calibrated foliations. An elementary family is obtained by applying
 the product construction to 3-dimensional taut foliations
and non-trivial symplectic manifolds.

Another important family of 2-calibrated foliations are cosymplectic foliations. Recall that they are given 
by a triple $(M^{2n+1},\alpha,\omega)$, 
$\alpha$ a closed 1-form and $\omega$ a closed 2-form such that $\alpha\wedge \omega^n$ is a volume form.
 An example of cosymplectic foliation  is a 2-calibrated foliation whose leaves
 are the fibers of a fibration over the circle;
the closed  1-form defining the foliation is the pullback of any volume form on the circle.
Each fiber is a closed symplectic manifold and the first 
return map associated to the kernel of the calibrating 2-form is a symplectomorphism.
We refer to such cosymplectic foliations as symplectic mapping tori. In fact, symplectic mapping tori are characterized 
 as cosymplectic foliations whose defining 1-form has rank one period lattice.
 This characterization implies that symplectic mapping tori are $C^0$-dense
in cosymplectic foliations. The reason is that the defining 1-form can be approximated
by closed 1-forms with rational periods.

Cosymplectic foliations appear naturally in symplectic geometry as follows: recall that a vector field $Y$ on a 
symplectic manifold  $(Z,\Omega)$ is called  symplectic if $L_Y\Omega=0$. If $Y$ is a symplectic vector field
 transverse to $\partial Z$, then its symplectic annihilator
\[\mathrm{Ann}(Y)^\Omega=\{v\in TZ\,|\,\Omega(Y,v)=0\}\]
is an integrable codimension one distribution. Since it contains the vector field $Y$, 
it induces a codimension one foliation $\mathcal{F}$
on $\partial M$. Let $\alpha:=i_Y\Omega$. It can be checked that $(\partial M,\alpha_{\partial M},\Omega_{\partial M})$
is a cosymplectic foliation. 

The previous construction leads to an analogy between cosymplectic foliations and contact structures. 
The reason is that on a symplectic manifold $(Z,\Omega)$ endowed with a vector field $Y$ transverse to the boundary and satisfying 
$L_Y\Omega=\Omega$, the restriction of $i_Y\Omega$ to $\partial M$ is a contact form. Following this analogy,
we define the Reeb vector  field $R$ of a cosymplectic foliation $(M,\alpha,\omega)$  to be
 the vector field characterized by the equations
$i_R\omega=0$, $i_R\alpha=1$. The foliation is invariant under the flow of the Reeb vector field.
 In fact, a cosymplectic foliation can be defined as a 2-calibrated foliation endowed with a vector 
field $R$ spanning the kernel of $\omega$ and whose flow preserves the foliation; we say that $R$ is a Reeb vector field.

A third family of 2-calibrated foliations are symplectic bundle foliations 
 \footnote{This family of 2-calibrated foliations was pointed out to the author by the referee.},
 which are defined as bundle foliations with fiber $S^1$
 over symplectic manifolds. There is a very rough way of associating symplectic bundle foliations 
to any bundle foliation  $\pi\colon M\rightarrow X$ with fiber $S^1$. The latter is characterized
by a conjugacy class of representations 
of $\pi_1(X,x)$ in $\mathrm{Diff}(S^1)$.  A result of Gompf (\cite{Go95}, theorem 0.1)
asserts that there exist closed symplectic manifolds (of dimension 4) whose fundamental group isomorphic to  $\pi_1(X,x)$.

\begin{example}\label{ex:3}
Let  $x_1,y_1,x_2,y_2,t$ be coordinates on  $\mathbb{R}^5$ and consider the canonical 2-form
$\omega_{\mathbb{R}^5}$. It descends to $\mathbb{T}^5=\mathbb{R}^5/\mathbb{Z}^5$ to a closed 2-form
$\omega_{\mathbb{T}^5}$. Let
$\mathcal{F}$ be any of the  foliations on $\mathbb{T}^5$ induced by
a constant 1-form $\alpha$ on $\mathbb{R}^5$ whose kernel is transverse to $\tfrac{\partial}{\partial t}$. Then
$(\mathbb{T}^5,\alpha,\omega_{\mathbb{T}^5})$ is a 2-calibrated foliation. Its leaves are 
all diffeomorphic to $\mathbb{R}^i\times \mathbb{T}^{4-i}$, where $i\in \{0,\dots,4\}$ depends on the slopes of the kernel
of the 1-form.

By construction  $(\mathbb{T}^5,\alpha,\omega_{\mathbb{T}^5})$ is both a
 cosymplectic foliation and a symplectic bundle foliation.
 It is a product (respectively a mapping torus) if and only if the leaves are diffeomorphic to 
 $\mathbb{R}^i\times \mathbb{T}^{4-i}$, $i\leq 2$ (respectively $\mathbb{T}^4$). 
\end{example}

Deciding which manifolds admit a 2-calibrated
foliation can be divided in several subproblems which in general are very hard.
A 2-calibrated foliation $(M,\mathcal{F},\omega)$ is the superposition of several
compatible structures. Firstly the foliation. Secondly the
2-form restricts to a closed non-degenerate foliated 2-form
$\omega_{\mathcal{F}}$. The pair $(\mathcal{F},\omega_\mathcal{F})$ defines a (regular) Poisson structure on $M$ and 
as such it is also defined by an appropriate bivector field $\Pi$. And thirdly the foliated
symplectic form  $\omega_{\mathcal{F}}$ admits a lift to a global closed 2-form
$\omega$.

Determining  which codimension one foliations are the
symplectic foliations of a Poisson structure is very complicated; there exist 
partial results which use h-principles
and only apply to open manifolds  \cite{Be01,Be02,FF10}. The existence of a closed lift of a
foliated
2-form $\omega_{\mathcal{F}}$ is controlled by three obstructions associated to the spectral
sequence which relates  basic cohomology,  leafwise cohomology and the cohomology of the total space \cite{EK} 
(see \cite{Al} for a treatment in the setting of Poisson geometry); if the foliation is defined by a closed 1-form, then 
the obstruction to the existence of a closed lift admits a simpler description (\cite{GMP10}, section 2.2).

We would like to regard a 2-calibrated foliation as a codimension one
regular Poisson manifold with a lift of $\omega_\mathcal{F}$ to a closed 
2-form $\omega$. We are not fully interested in the  2-form $\omega$, as the
following definition reflects.

\begin{definition}\label{def:equiv} Let $(M_j,\mathcal{F}_j,\omega_j)$, $j=1,2$, be 
2-calibrated foliations. They are said to be
equivalent if there exists a diffeomorphism $\phi\colon
M_1\rightarrow M_2$ such  that
\begin{itemize}
\item $\phi$ is a Poisson morphism or equivalence (it preserves the
foliations together with the leafwise 2-forms);
 \item $[\phi^*\omega_2]=[\omega_1]\in H^2(M_1;\mathbb{R})$ and $\phi$ preserves the co-orientations.
\end{itemize}
\end{definition}

For symplectic mapping tori an equivalence 
is just a Poisson diffeomorphism preserving co-orientations.
Alternatively, equivalent symplectic mapping tori are those with the same symplectic
leaf an isotopic first return maps (the isotopy through symplectomorphisms).

As we shall see in the following sections, the notion of equivalence is the right one 
to remove the dependence on choices in our surgeries.

\section{Normal Connected Sum}\label{sec:nsum}
In the previous section we saw that deciding whether a manifold supports a
2-calibrated
foliation is very complicated. It is thus natural to look for procedures to
 build new 2-calibrated foliations out of given ones. In this section we introduce the normal connected sum 
of integral 2-calibrated foliations, and we use it to give examples of 2-calibrated foliations which do not
belong to either of the three elementary families, hence proving proposition \ref{pro:newex}.

Symplectic normal connected sum is a surgery construction in
which two symplectic manifolds are glued along two copies of the
same codimension two symplectic submanifold, which enters in the
manifolds with opposite normal bundles (\cite{Go95}, theorem 1.3).
A parametric version of this surgery gives rise to an analogous construction for regular Poisson
manifolds (\cite{IM03}, theorem 1). We propose the following extension to integral 2-calibrated foliations.

\begin{theorem}\label{thm:calibnorsum} Let $(M^{2n+1}_j,\mathcal{F}_j,\omega_j)$, $j=1,2$, be
 integral
2-calibrated foliations. Let $(N^{2n-1},\mathcal{F}_N,\omega_N)$ be a
2-calibrated foliation which is a symplectic mapping torus. Assume that we have maps $l_j\colon
N\hookrightarrow M_j$, $j=1,2$, embedding $N$ as a 2-calibrated submanifold of
$M_j$ (definition \ref{def:submnfd}),  such that the following properties hold:
\begin{enumerate}
\item The  2-calibrated foliations induced by the embeddings are equivalent to
the given one $(N,\mathcal{F}_N,\omega_N)$ (definition \ref{def:equiv}).
\item  The normal bundles of $l_j(N)\subset M_j$, $j=1,2$, are trivial.
\item The fiber of $N\rightarrow S^1$ is simply connected.
\end{enumerate}

Then there exist gluing maps $\psi$ such the Poisson
structure $\Pi$ on ${M_1}\#_\psi M_2$  characterized by matching on $M_j\backslash l_j(N)$ 
the Poisson structures  $\Pi_j$ 
associated to $(M_j,\mathcal{F}_j,\omega_j)$, $j=1,2$, admits a lift to a 2-calibrated structure.

\end{theorem}
\begin{proof}
By assumptions 1 and 2 Poisson surgery produces a Poisson structure $\Pi$ on
${M_1}{\#_\psi}M_2$ \cite{IM03}. Very briefly, there is a gluing map $\psi$ identifying
 $A_1\rightarrow A_2$ annular neighborhoods of
$l_1(N)$ and $l_2(N)$ (by this we mean tubular neighborhoods from which we remove $l_j(N)$, $j=1,2$) defined as follows:
by assumption 2 the normal bundles are trivial and by
Darboux-Weinstein theorem with parameters the (smooth) leaf space of $N$ (\cite{MS98}, chapter 3), there exist trivializations
in which $\Pi_j$, $j=1,2$, split. One factor is the leafwise symplectic form on $l_j(N)$ and the other one is the 
standard symplectic form $dx\wedge dy$ on the normal disk with coordinates $x,y$.
On each normal disk $\psi$ 
is the unique rotationally independent symplectomorphism of the punctured disk of radius $\delta>0$
which reverses the orientation of the radii.

Let $(\mathcal{F},\omega_\mathcal{F})$ denote the foliation and leafwise symplectic form associated
 to $\Pi$. If there is a lift of $\omega_\mathcal{F}$
 to an integral closed 2-form $\omega$, then there must be a Hermitian line bundle $L$
and a compatible connection $\nabla$ such that 
\[-2\pi i\omega=F_\nabla,\]
where $F_\nabla$ is the curvature of the connection.

Because $w_j$, $j=1,2$, represent integral cohomology classes there exist $(L_j,\nabla_j)\rightarrow M_j$ Hermitian 
line bundles with compatible connections such that
\begin{equation}\label{eq:lcurv}
-2\pi i\omega_j=F_{\nabla_j}.
\end{equation}
 We look for a
 lift of $\psi$ to a bundle isomorphism
\[\Psi\colon {L_1}_{\mid A_1}\rightarrow {L_2}_{\mid A_2},\] to define a (Hermitian) line bundle
 \[L:=L_1{\#_\Psi} L_2\rightarrow  {M_1}\#_\psi M_2.\]
Let $c_j$, $j=1,2$, denote the Chern classes of ${L_j}_{\mid A_j}$, 
which are integral lifts of the restrictions of $w_j$ to $A_j$. 
An isomorphism lifting $\psi$ exists if and only if
\begin{equation}\label{eq:integralcoh}
 \psi^*c_2=c_1 \in H^2(A_1;\mathbb{Z}).
\end{equation}
 Because the fiber of $N\rightarrow S^1$ is simply connected, 
the Wang sequence for the mapping torus $A_1\rightarrow S^1$
implies that $H^2(A_1;\mathbb{Z})$ is torsion free. Therefore   
equation (\ref{eq:integralcoh}) is equivalent to 
\begin{equation}\label{eq:realcoh}
 [\psi^*{w_2}_{\mid A_2}]=[{w_1}_{\mid A_1}]\in H^2(A_1;\mathbb{R}).
\end{equation}
Because $w_j$, $j=1,2$,  extend to $A_j\cup l_j$ and the cohomology of the tubular neighborhoods
is concentrated in 
$l_j(N)$, equation (\ref{eq:realcoh}) is equivalent to 
\[ [l_2^*w_2]=[l_1^*w_1] \in H^2(N;\mathbb{R}),\]
which holds true because by assumption 1 the 2-calibrations induced by $l_1$ and $l_2$ on $N$ are equivalent.

Therefore we obtain $L\rightarrow {M_1}{\#_\psi}M_2$ a Hermitian line
bundle with two not everywhere defined compatible connections
$\nabla_1,\nabla_2$, overlapping on $A_1\subset {M_1}{\#_\psi}M_2$. Remark that by equation (\ref{eq:lcurv}) 
the leafwise curvatures match on $A_1$. We are going to use the assumptions to modify $\nabla_1$ and $\nabla_2$ (the latter 
away from $l_2(N)$),
so that we obtain the leafwise equality of connections on $A_1$. Then a convex  combination
of both connections associated to a partition of the unity subordinated to $M_j\backslash l_j(N)$, $j=1,2$,
  is a connection on  ${M_1}{\#_\psi}M_2$ whose leafwise curvature is $-2\pi i\omega_\mathcal{F}$.

The difference
\begin{equation}\label{eq:leafconn}
l_1^*\nabla_1-l_2^*\nabla_2
\end{equation} is a leafwise closed 1-form on $N$ (recall that $N$ is a mapping torus and therefore all
 leaves are compact). By assumption 3 it is leafwise
exact and therefore we can modify say $\nabla_2$,
 by adding a smooth leafwise primitive function
 so the 1-form in equation
(\ref{eq:leafconn}) is leafwise vanishing.

 Triviality of the normal bundles implies
the existence of normal
forms for the leafwise connections on tubular neighborhoods of $l_j(N)$, $j=1,2$,
 which only depend on the restrictions of the leafwise
connections to  $l_j(N)$; the normal forms amount to fixing a primitive 1-form for $dx\wedge dy$.
 The connections can be
assumed to coincide with the normal forms. Finally the difference
$\nabla_1-\psi^*\nabla_2$
is not still leafwise vanishing; on each normal annulus it is the differential of an (explicit) function, and what we do 
is modifying accordingly $\nabla_2$ on $M_2\backslash l_2(N)$.

As for dependence of the construction on choices, remark that the choice of isotopy 
classes of trivializations of the normal bundles (the framings),
may affect the diffeomorphism class of $M_1\#_\psi M_2 $. For fixed isotopy classes  of trivializations of the normal bundles,
the underlying Poisson structure is unique up to Poisson diffeomorphism. The reason is that 
the leafwise symplectic form is unique up to isotopy supported near $N$. This follows from an
elementary argument
which is going to be used several times: because the leaves of $N$
have no first cohomology group the local path of symplectomorphisms provided by Moser's argument is Hamiltonian (\cite{MS98},
 chapter 3). The choice of primitive Hamiltonian function can be done coherently for all leaves of $N$. 
By extending the corresponding
function to a global one supported near $N$, we construct a path of transformations connecting both Poisson structures.
Also, if we  fix a isotopy class of lifts $\Psi$, the 2-calibrated structure provided by the normal
connected sum is unique up to equivalence. This is because the cohomology class of the calibrating 2-form is the
 image in real cohomology of the first Chern class
of the bundle $L$, which is fixed by the choice of isotopy class of lifts. 
\end{proof}

\begin{remark} The hypothesis needed to define normal connected sum of regular Poisson manifolds  are much
weaker than the requirements in theorem \ref{thm:calibnorsum}. In particular the normal bundles $l_j(N)$, $j=1,2$,
are not required to be trivial, just opposite. Triviality of the normal bundles
 is necessary if we want to produce an integral 2-calibrated
foliation extending the given Poisson structures $\Pi_j$ on $M_j\backslash l_j(N)$, $j=1,2$. The reason is 
that already in the symplectic setting, having non-trivial normal bundle gives rise to choices 
in the construction which result into symplectic forms with different volume; 
this is a well known issue that appears when blowing up
 symplectic submanifolds (\cite{MS98}, chapter 7). 

Perhaps the assumptions in theorem \ref{thm:calibnorsum} can be weakened if we just 
require the existence of a 2-calibration on the normal 
connected sum.
\end{remark}

The normal connected sum can be applied to construct integral 2-calibrated
foliations, that use as building blocks 2-calibrated foliations
which are products and symplectic mapping tori,  but which are neither products, nor cosymplectic foliations
nor symplectic bundle foliations.

\begin{proof}[Proof of proposition \ref{pro:newex}]
Let $(P^4,\Omega)$ be an integral symplectic 4-manifold which contains
a symplectic sphere $S^2$ with trivial normal bundle; let $A\in\mathbb{Z}$ be the
induced area form on the sphere. Let $\varphi\in
\text{Symp}(P,\Omega)$  such that $\varphi_{\mid S^2}=\text{Id}$; for example $\varphi$ can be the identity. We define
$(M_1,\mathcal{F}_1,\omega_1)$ to be the symplectic mapping torus associated to $\varphi$.

Let $(M_2,\mathcal{F}_2,\omega_2)$ be the product  2-calibrated foliation
with factors any taut foliation $(Y^3,\mathcal{F}^3,\sigma)$ and the sphere $(S^2,A)$;
 via a small perturbation and a rescaling of $\sigma$,
we may take $\omega_2$ to be integral. Let $C$ be a
fixed transverse cycle for $(Y^3,\mathcal{F}^3,\sigma)$ and $\theta\colon
S^1\rightarrow C$ any fixed positive parametrization with respect to  the
co-orientation.

Let $N^3$ be the result of applying the mapping torus construction to
$\mathrm{Id}\in \mathrm{Symp}(S^2,A)$ ($N\cong S^1\times S^2$). Since
$\varphi_{\mid
S^2}=\text{Id}$, there is an obvious embedding $l_1\colon
N\hookrightarrow M_1$.  The embedding $l_2$ is the product map
$\theta\times \text{Id}\colon N\hookrightarrow M_2$.

By construction the embeddings  fulfill the hypothesis of theorem
\ref{thm:calibnorsum}, so we obtain a 2-calibrated foliation
$({M_1}{\#_\psi}M_2,\mathcal{F},\omega)$.

We impose the following additional constraints on the summands to make sure that
 $({M_1}{\#_\psi}M_2,\mathcal{F},\omega)$ does not belong to the three
basic families:
\begin{itemize}
 \item  $(Y^3,\mathcal{F}^3)$ contains compact and non-compact leaves.
\item There is a compact leaf $\Sigma$ of $(Y^3,\mathcal{F}^3)$ which intersects $C$ in exactly one point, and $(P^4,\Omega)$
is an odd Hirzebruch surface  (\cite{MS98}, chapter 4).
\item The genus of $\Sigma$ is greater than one, and $\pi_1(Y)$ is not isomorphic to $\pi_1(S^1\times\Sigma)$.
\end{itemize}
 
Because $l_2(N)$ intersects each leaf of $(M_2,\mathcal{F}_2)$ in a unique connected component, 
there is a one to one correspondence between  leaves of
$(Y^3,\mathcal{F}^3)$ and  leaves of $({M_1}{\#_\psi}M_2,\mathcal{F})$. This correspondence sends a leaf $F$ of
$(Y^3,\mathcal{F}^3)$ to the leaf which contains $(F\times S^2)\backslash (l_2(N)\cap (F\times S^2))$. 
Because the leaves of $(M_2,\mathcal{F}_2)$ are compact,
the  correspondence sends compact leaves to compact leaves and
non-compact leaves to non-compact leaves. Since $(Y^3,\mathcal{F}^3)$ contains compact and non-compact leaves
so  does $({M_1}{\#_\psi}M_2,\mathcal{F},\omega)$,  and hence it has non-trivial holonomy.
Consequently, $({M_1}{\#_\psi}M_2,\mathcal{F},\omega)$ cannot be a cosymplectic foliation.

Let $\Sigma$ be a compact leaf of $\mathcal{F}^3$ which intersects $C$ in one point. The correspondence between
 leaves described in the previous paragraph
sends $\Sigma$ to a compact leaf $F_\Sigma$, which is the symplectic normal connected sum of the odd Hirzebruch surface and
$(\Sigma\times S^2,p_1^*\omega_{\mid \Sigma}+ p_2^*A)$ 
along a symplectic sphere with trivial normal bundle. At the differentiable level $F_\Sigma$ 
is the normal connected sum of the trivial $S^2$-fibration
 over $\Sigma$ and the twisted
$S^2$-fibration over $S^2$, and hence it is the twisted $S^2$-fibration over $\Sigma$ (the fibers of our fibrations have a coherent 
orientation, since they are symplectic). If $F_\Sigma$ is diffeomorphic to a product of surfaces 
then we can only have $F_\Sigma\cong S^2\times \Sigma$; otherwise we could not have isomorphic fundamental groups. But
then $F_\Sigma$ would admit two different $S^2$-fibration structures, and this is in contradiction with  
\cite{Me84}. Therefore $({M_1}{\#_\psi}M_2,\mathcal{F},\omega)$ cannot be a product.

If the normal connected sum is a symplectic bundle foliation $\pi\colon {M_1}{\#_\psi}M_2\rightarrow X$, then 
$F_\Sigma$ is a covering space of $X$. Because the fundamental group of $F_\Sigma$ is the fundamental group of $\Sigma$,
 our assumption
on the genus of $\Sigma$ implies that the covering must be trivial. Therefore $\pi$ sends $F_\Sigma$ diffeomorphically onto $X$.
This also implies
that the principal $S^1$-bundle has a section, so ${M_1}{\#_\psi}M_2$ is the trivial bundle $S^1\times F_\Sigma$. Hence
$\pi_1({M_1}{\#_\psi}M_2)$ is diffeomorphic to $\pi_1(S^1\times \Sigma)$. But applying
 Seifert-Van Kampen theorem to the open subsets $M_1\backslash l_1(N)$, $M_2\backslash l_2(N)$ gives
that $\pi_1({M_1}{\#_\psi}M_2)$ is diffeomorphic to $\pi_1(Y)$, and this contradicts the assumption on $\pi_1(Y)$.

\end{proof}

\section{Generalized Dehn Surgery}\label{sec:dsurg}
In this section we introduce our second surgery,  generalized Dehn surgery. We give a first
definition which is the most natural one from the viewpoint of foliation theory. We present a second
approach via handle attaching along Lagrangian spheres; this is a very natural definition having into account
the description of Legendrian surgeries in contact geometry (\cite{We91}, section Elementary Cobordisms).
 We prove the equivalence of both constructions 
in theorem \ref{thm:equiv}.

Generalized Dehn surgery is done, unlike  normal
connected sum, along a submanifold inside one of the leaves.
Let $(M,\mathcal{F},\omega)$ be a 2-calibrated foliation. We orient $M$ so that
a positive transverse vector  followed
by a positive basis of the leaf with respect to the volume form
 $\omega_{\mathcal{F}}^n$, gives a positive basis.

Let $T:=T^*S^n$ and  $d\alpha_{\mathrm{can}}$ its canonical symplectic structure. Let \[\tau\colon T\rightarrow T\]
be a generalized Dehn twist. Recall that these are certain compactly
supported symplectomorphisms of $(T,d\alpha_{\mathrm{can}})$ which induce the
antipodal map on the zero section.  Let $T(\lambda)$ be the subset of cotangent vectors of length
$\leq \lambda$ with respect to  the round metric. Generalized Dehn twists can be chosen
to be supported in the interior of $T(\lambda)$ for any fixed $\lambda$, and
any two with such property are isotopic  in
$\mathrm{Symp}^{\mathrm{comp}}(T(\lambda),d\alpha_{\mathrm{can}})$, the group of compactly supported symplectomorphisms
 (\cite{Se03}, lemma 1.10 in section 1.2). They are symplectic generalizations of Dehn
twists on $T^*S^1$.

A parametrized Lagrangian sphere $L\subset (M,\mathcal{F},\omega)$  is a submanifold of a 
leaf $F_L$ such that  $\omega_{ L}\equiv 0$,
 together with a parametrization $l\colon
S^{n}\rightarrow L$. 
By a theorem of Weinstein (\cite{MS98}, chapter 3)
there exists $U$ a compact
neighborhood of   $L$ inside $F_L$
 and $\lambda>0$, such that   $l^{-1}\colon
L\rightarrow S^n$ extends to a
symplectomorphism
\[
 \varphi\colon(U,\omega_{\mathcal{F}})\rightarrow (T(\lambda),d\alpha_{\mathrm{can}}).
 \]
Let us assume that if $n=1$ the loop $L$ has trivial holonomy; if $n>1$ the absence of holonomy is a consequence of Reeb's theorem. 
In a neighborhood of $L$ the foliation is a product. We let $R$
be a local positive Reeb vector field and we let $\varPhi^R_t$ denote its time $t$ flow, which
by definition preserves $\mathcal{F}$. Let $\epsilon>0$ small enough
so that
\begin{eqnarray*}
\varPhi^R\colon [-\epsilon,\epsilon]\times U &\longrightarrow & M\\
(t,x)&\longmapsto &\varPhi^R_t(x)
\end{eqnarray*}
is an embedding. We introduce the following notation:

\begin{eqnarray}\label{eq:not}\nonumber
U(\epsilon):=\varPhi^R( [-\epsilon,\epsilon]\times U), & U_t:=\varPhi^R_t(U),\\
U^+(\epsilon):=\varPhi^R( [0,\epsilon]\times U),&
U^-(\epsilon):=\varPhi^R( [-\epsilon,0]\times U).
\end{eqnarray}
The result of cutting $U(\epsilon)$ along $U$ is the manifold $U^-(\epsilon)\coprod U^+(\epsilon)$ whose boundary
contains  $U^-=U\times \{0\}\subset U^-(\epsilon)$, $U^+=U\times \{0\}\subset U^+(\epsilon)$. 

\begin{definition}\label{def:gdehn} Let $L\subset (M,\mathcal{F},\omega)$ be a parametrized Lagrangian sphere.
If $n=1$ assume that
 $L$ is a loop with trivial holonomy. Generalized Dehn surgery  along $L$ is defined
 by cutting $M$  along $U$ as above and then gluing back via the
composition

\begin{equation}\label{eq:Dehn} \chi\colon
(U^-,\omega_{\mathcal{F}})\overset{\varphi}{\rightarrow}(T(\lambda),d\alpha_{\mathrm{can}}
)\overset{\tau}
{\rightarrow}
(T(\lambda),d\alpha_{\mathrm{can}})\overset{\varphi^{-1}}{\rightarrow}(U^+,\omega_{
\mathcal{F}}),
\end{equation}
where $\tau$ is any  choice of generalized Dehn twist
supported in the interior of $T(\lambda)$ and we use the canonical identifications of $U^-,U^+$ with $U$.

We denote the resulting foliated manifold by $(M^L,\mathcal{F}^L)$.
\end{definition}

\begin{proposition}\label{pro:gDehn} The foliation $(M^L,\mathcal{F}^L)$ admits calibrations $\omega^L$. 
 If $n>1$ then
\begin{enumerate}
\item $(M^L,\mathcal{F}^L,\omega^L)$ is unique up to equivalence;
\item $[\omega]$ is integral if and only if $[\omega^L]$ is
 integral;
 \item  $\pi_i(M^L)\cong\pi_i(M)$ and $H_i(M^L;\mathbb{Z})\cong
 H_i(M;\mathbb{Z})$, $0\leq i\leq n-1$.
 \end{enumerate}
\end{proposition}

\begin{proof} 
We restrict our attention to $U(\epsilon)$. After cutting  $U(\epsilon)$ along $U$ and gluing back
using the identification
$\chi$ in equation (\ref{eq:Dehn}), we obtain
\[U^L(\epsilon):=U^-(\epsilon){\#_\chi} U^+(\epsilon)\subset M^L.\]
Since the flow of $R$ preserves both $\omega$ and the foliation, the restriction of $\omega$
to $U^-(\epsilon)$ and $U^+(\epsilon)$ defines closed 2-forms $\omega^-$
and $\omega^+$ independent of the coordinate $t$. When we glue $U^-$ to
$U^+$ using $\chi$, being this map a symplectomorphism the 2-forms
$\omega^-$ and $\omega^+$ induce on $U^L(\epsilon)$ a 2-form
$\omega^L_\epsilon$. Then
\begin{equation*}
 \omega^L:=
 \begin{cases}
 \omega & \text{in}\;\; M^L\backslash U^L(\epsilon),\\
 \omega^L_\epsilon & \text{in}\;\; U^L(\epsilon)
 \end{cases}
 \end{equation*}
 is the desired closed 2-form.

The 2-calibrated structure we obtain is unique up to equivalence. Firstly different identifications
$\varphi\colon(U,\omega_{\mathcal{F}})\rightarrow (T(\lambda),d\alpha_{\mathrm{can}})$ are related by a global Poisson 
diffeomorphism. The reason is the same as in the proof of the uniqueness statement of theorem \ref{thm:calibnorsum}:
 $S^n$, $n>1$, is simply connected. Secondly
  generalized Dehn twists are
symplectically isotopic by an isotopy
supported in a  neighborhood of the sphere. Thirdly changing the Reeb vector field
amounts to a change of variable in the coordinate $t$, and this does not modify the construction. 

The calibration is a real cohomology class determined by its values on closed 2-chains (which by a theorem of Thom 
are always homologous to embedded surfaces). If $n>1$ the 2-chains can be homotoped to avoid the neighborhood $U(\epsilon)$
of the Lagrangian sphere $L$, where $\omega^L$ coincides with $\omega$. Hence the integrality of the 2-calibrated
foliation is unaffected by the surgery.  

The same general position arguments imply that maps from CW complexes of dimension less or equal than $n$ can be homotoped to miss
$U(\epsilon)$. Therefore homology and homotopy groups up to dimension $n-1$ are unaffected by the surgery.

\end{proof}
\begin{remark}\label{rem:framed} A ``framed" Lagrangian $n$-sphere \cite{Se00} is a
parametrized $n$-sphere up to isotopy
and the action of $O(n+1)$. Generalized Dehn twists associated
to two parametrizations defining the same ``framed" Lagrangian $n$-sphere are
isotopic, the isotopy by symplectomorphisms supported in a compact
neighborhood of the Lagrangian sphere (remark 5.1 in \cite{Se00} or paragraph after lemma 1.10 in \cite{Se03}).
Therefore generalized Dehn surgery is well defined for ``framed"
Lagrangian spheres.
\end{remark}

\begin{remark}\label{rem:slide}  The flow of the local Reeb vector field $R$ can be used to displace the
Lagrangian sphere $L$ to a new
 Lagrangian sphere $L'$ inside a nearby leaf. It follows that $(M^L,\mathcal{F}^L,\omega^L)$ and
$(M^{L'},\mathcal{F}^{L'},\omega^{L'})$ are equivalent.

If we use instead of $\tau$ its inverse, we
 get a  2-calibrated foliation $(M^{L^-},\mathcal{F}^{L^-},\omega^{L^-})$ referred
 to as negative generalized Dehn surgery along $L$; negative
 generalized Dehn surgery is generalized Dehn surgery for the opposite
co-orientation.

  Generalized Dehn
  surgery along $L$ and negative generalized Dehn surgery
   along $L$  are inverse of each other.
\end{remark}

\subsection{Lagrangian surgery}\label{ssec:lagsurg} Let $L\subset (M,\mathcal{F},\omega)$ be a parametrized Lagrangian sphere,
and let $\nu(L)$ and $\nu_\mathcal{F}(L)$ denote respectively a  tubular neighborhood of $L$ and a 
tubular neighborhood of $L$ inside the leaf
 containing $L$. The parametrized Lagrangian sphere $L$
carries a canonical framing $\mu_L$: because $L$ is Lagrangian $\nu_\mathcal{F}(L)\cong
T^*L$ and we deduce
\begin{equation}\label{eq:framings}
\nu(L)\cong \nu_\mathcal{F}(L)\oplus \underline{\mathbb{R}}\cong T^*S^n\oplus
\underline{\mathbb{R}}\cong {\underline{\mathbb{R}}^{n+1}}_{\mid S^{n}},
\end{equation}
where in the last
isomorphism in (\ref{eq:framings}) a positive no-where vanishing section of
$\underline{\mathbb{R}}_{\mid S^{n}}$ is sent to the outward normal unit vector field.
Therefore $L$ determines up to diffeomorphism an elementary cobordism $Z$, which amounts to attaching a $(n+1)$-handle
to the parametrized sphere $L$ with  framing $\mu_L$ (\cite{GST00}, chapter 4). The boundary of the cobordism is
$\partial Z=M\coprod M^{\mu_L}$.

This subsection  addresses the construction of a 2-calibrated foliation  
$(M^{\mu_L},\mathcal{F}^{\mu_L},\omega^{\mu_L})$  which extends $(M,\mathcal{F},\omega)$ on the 
complement of a neighborhood of $L$ (the complement understood as a subset of both $M$ and $M^{\mu_L}$). We do it by using the
relation between symplectic manifolds and cosymplectic foliations presented in section \ref{sec:def}:
  we have to endow the cobordism $Z$ with
a symplectic form $\Omega$ -at least in a neighborhood of the $(n+1)$-handle- and a symplectic
 vector field $Y$ transverse to the
boundary. This produces automatically a cosymplectic foliation on $\partial Z$, and that is how we obtain
$(M^{\mu_l},\mathcal{F}^{\mu_L},\omega^{\mu_L})$.
Remark that our strategy is the same one used in
contact geometry to show that surgeries along Legendrian spheres give
rise to new contact manifolds (\cite{We91}, paragraph 3 in page 242).

The elementary cobordism $Z$ is the result of gluing a $(n+1)$-handle to the trivial
cobordism $P_1:=M\times  [-\varepsilon,\varepsilon]$. We have to define symplectic
structures and symplectic vector fields transverse to the boundary
 on both the trivial cobordism and the $(n+1)$-handle, in 
a way that is compatible with the gluing.

We start with the trivial cobordism $P_1$:
by the coisotropic embedding \cite{Go83} there is a unique choice
of symplectic structure on $P_1$ which extends the given closed 2-form $\omega$  on $M\times \{0\}$.
 We now give  a specific normal
form for it  which is convenient for the purpose of describing a compatible gluing with the $(n+1)$-handle:
let us denote $H_1:=\nu(L)$. Since
the gluing between the trivial cobordism and the $(n+1)$-handle occurs near $\nu(L)$,  we can 
assume without loss of generality that $P_1=H_1\times  [-\varepsilon,\varepsilon]$. Let
 $(\mathcal{F}_1,\omega_1)$ denote the restriction of $(\mathcal{F},\omega)$ to $H_1$. We select
$R_1$ a positive Reeb vector field on $H_1$  with  dual (closed) defining 1-form 
$\alpha_1$ ($i_R\alpha_1=1$, $\mathrm{ker}\alpha_1=\mathcal{F}_1$). 
 We let $v$ be the coordinate on the interval $[-\varepsilon,\varepsilon]$, and we extend $\alpha_1$
 and $\omega_1$ to $H_1\times  [-\varepsilon,\varepsilon]$  independently of $v$.

We define on $P_1$
\[\Omega_1:=\omega_1+d(v\alpha_1),\]
which is a symplectic form provided $\varepsilon$ is small enough.

As symplectic vector field on $(P_1,\Omega_1)$ we take
$Y_1:=\tfrac{\partial}{\partial v}$, which is  transverse to $H\times \{-\varepsilon\}$ and $H\times \{\varepsilon\}$.

We let $P_2$ denote the $(n+1)$-handle. Before defining the symplectic form $\Omega_2$ and a symplectic vector field $Y_2$ on
$(P_2,\Omega_2)$,  we address
the problem of gluing symplectic cobordisms.

\begin{lemma}[\cite{Go83}, Extension theorem] \label{lem:cosympglue}
Let $(P_j,\Omega_j)$, $j=1,2$, be symplectic manifolds, $H_j\subset P_j$ hypersurfaces
and $Y_j$ symplectic vector
fields  transverse to them, so
that we have product structures $H_j\times
[-\varepsilon,\varepsilon]$. Define $\omega_j={\Omega_j}_{\mid H_j}$,
$\alpha_j={i_{Y_j}\Omega_j}_{\mid H_j}$ and $\mathcal{F}_j$ the foliation integrating $\mathrm{ker}\alpha_j$, $j=1,2$.
  Suppose that $\phi\colon H_1\rightarrow H_2$ is
   a diffeomorphism such that $\phi^*\omega_2=\omega_1$ and
$\phi^*\alpha_2=\alpha_1$ (and therefore  $\phi^*\mathcal{F}_2=\mathcal{F}_1$). Then
  \[ \phi\times \text{Id}\colon (H_1\times
[-\varepsilon,\varepsilon],\Omega_1)\rightarrow
(H_2\times [-\varepsilon,\varepsilon],\Omega_2)\] is a
symplectomorphism (obviously compatible with the symplectic vector fields).
\end{lemma}
Lemma \ref{lem:cosympglue} is the analog of proposition 4.2 in \cite{We91}.

In our specific situation of gluing near Lagrangian spheres, the amount of information needed 
 to describe $\phi$ as in 
lemma \ref{lem:cosympglue} is much smaller.
\begin{corollary}\label{cor:cosympglue} Let $(P_j,\Omega_j,H_j,Y_j)$, $j=1,2$, be as
in lemma \ref{lem:cosympglue} and assume further that $L_j\subset H_j$ are Lagrangian
spheres  and  $P_j$  small tubular neighborhoods of $L_j$.

Let $\theta\colon L_1\rightarrow L_2$ be a diffeomorphism.
Then $\theta$ extends to an isomorphism of tuples
 \[(P_1,\Omega_1,H_1,Y_1)\rightarrow (P_2,\Omega_2,H_2,Y_2).\]
\end{corollary}
\begin{proof}  The symplectic vector fields give rise by contraction
to closed 1-forms defining the foliations, and therefore to Reeb vector fields.
We extend $\theta$ to a symplectomorphism of neighborhoods of the spheres inside their
 leaves,  and we further extend it to  $\phi\colon
(H_1,\alpha_1,\omega_1)\rightarrow (H_2,\alpha_2,\omega_2)$ by declaring it to be equivariant with respect to the Reeb
flows. By construction $\phi$ is in the hypothesis of
lemma \ref{lem:cosympglue}.

Notice that the only choice is the identification of the symplectic 
neighborhoods of $L_j$, $j=1,2$, inside their respective leaves. 
\end{proof}

\subsubsection{The choice of symplectic form and symplectic vector field on the
$(n+1)$-handle.}\label{sssec:shape}

Let $W$ be a neighborhood of $0\in \mathbb{C}^{n+1}$. This neighborhood will  contain our $(n+1)$-handle $P_2$.  

Let us consider the complex Morse function
 \begin{eqnarray*}
 h\colon \mathbb{C}^{n+1} &\longrightarrow &\mathbb{C} \\
 (z_1,\dots,z_{n+1})&\longmapsto &z_1^2+\cdots +z_{n+1}^2.
\end{eqnarray*}
We take $\Omega_2\in \Omega^2(W)$ to be any symplectic form of type $(1,1)$
at the origin with respect to the standard complex structure of $\mathbb{C}^{n+1}$,
and $Y_2$ to be the Hamiltonian vector field of $-\mathrm{Im}h$.

Let us explain the reason behind the choice of $(\Omega_2,Y_2)$. In
 the construction of the symplectic $(n+1)$-handle we have to reconcile 
several aspects: 

The  data $(P_2,\Omega_2,Y_2)$ has to determine the standard $(n+1)$-handle:  
if $\Omega_2=\Omega_{\mathbb{R}^{2n+2}}$ then $Y_2$ is 
the gradient flow of $-\mathrm{Re}h$ with respect to the Euclidean metric, whose dynamics determine the standard $(n+1)$-handle.  
In  lemma \ref{lem:stb} we are going to prove that for  $\Omega_2$ of type $(1,1)$ at the origin, the Hamiltonian vector field $Y_2$
has a  hyperbolic singularity at $0\in \mathbb{C}^{n+1}$. Therefore the flow of $Y_2$ has both 
the right dynamical behavior to construct a standard
$(n+1)$-handle about $0\in \mathbb{C}^{n+1}$ and the right symplectic behavior.

The second aspect is that we want to define Lagrangian surgery along $L$ so that
 it becomes equivalent to generalized Dehn surgery. Generalized Dehn twists 
appear in our current setting as follows: the origin $0\in \mathbb{C}^{n+1}$ is an isolated critical point for $h$. 
Let $h_z$ denote the fiber $h^{-1}(z)\cap W$, $z\in \mathbb{C}$, and let $\Omega$ be any
closed 2-form on  $W$ for which the fibers $h_z$ are symplectic. The
 annihilator with respect to $\Omega$ of the tangent space to the fibers is an 
Ehresmann connection for
$h\colon W\backslash \{0\} \rightarrow \mathbb{C}$. Parallel transport over a path 
not containing the critical value $0\in \mathbb{C}$,
defines a symplectomorphism from the regular
fiber over the starting point to the regular fiber over the  ending point. 
Seidel proves (\cite{Se03}, lemma 1.10 in section 1.2)
that for certain choice of closed 
2-form  $\Omega_\tau$ which is Kahler near the origin and for all $r\in \mathbb{R}^{>0}\subset \mathbb{C}$, 
 parallel transport of the fiber $h_r$
 over the boundary of the disk $\overline{D}(r)\subset \mathbb{C}$  counterclockwise, is
conjugated to a generalized Dehn twist supported in a given $T(\lambda)$. An argument using Taylor expansions shows
that for symplectic forms of type $(1,1)$ at the origin the fibers  $h_z$ are symplectic near the origin, 
and therefore there is an associated symplectic
parallel transport with respect to $\Omega_2$. Besides, symplectic parallel transport  with respect to $\Omega_2$
can be connected to symplectic parallel transport with respect to $\Omega_\tau$. The upshot
is that  symplectic parallel transport over $\overline{D}(r)\subset \mathbb{C}$ counterclockwise 
with respect to $\Omega_2$ can be isotoped to a generalized Dehn twist,
which is the property we need to prove the equivalence of generalized Dehn surgery and Lagrangian surgery.

The third aspect is that we need a flexible choice of symplectic form $\Omega_2$ on the $(n+1)$-handle,
 so the cobordisms naturally associated
to Lefschetz pencil structures to be described in subsection \ref{ssec:surgfiber}, can be identified with Lagrangian surgery.

In the next lemma we collect some useful properties of parallel transport with respect to forms of type $(1,1)$ at the origin:

\begin{lemma}\label{lem:stb} Let $\Omega\in \Omega^2(W)$ be a symplectic form
of type $(1,1)$ at the origin. Let $Y\in
\mathfrak{X}(W)$ be the Hamiltonian vector field of $-\mathrm{Im}h$ with respect to 
$\Omega$. Then the following holds:
\begin{enumerate}
\item $Y$ is a section of $\mathrm{Ann}(Y)^{\Omega}$ which vanishes at $0\in \mathbb{C}^{n+1}$.
\item $h_*Y(p)$ is a strictly negative multiple of $\tfrac{\partial}{\partial x}$, where $p\in W\backslash\{0\}$, $z=(x,y)$.
\item $Y$ has a non-degenerate singularity at the origin with $n+1$
positive eigenvalues
and $n+1$ negative eigenvalues.
\item  For each $r\in \mathbb{R}\backslash \{0\}$ we have Lagrangian spheres
$\Sigma_{r}\subset h_r$ characterized as the set of
points  contracting into the critical point by the parallel
transport over the segment $[0,r]$; the spheres come with a
parametrization up to isotopy and the action of $O(n+1)$ (they are ``framed").
More generally, for each $z$ and $\gamma$
an embedded curve joining $z$ and the origin, the points in $h_{\gamma(0)}$
sent to the origin by parallel transport over $\gamma$ are
a Lagrangian sphere $\Sigma_{\gamma(0)}$. Their construction depends smoothly on $\Omega$ and
$\gamma$.
\item For any embedded curve $\gamma$ through the origin  parallel transport
\[\rho_\gamma\colon (h_{\gamma(0)}\backslash \Sigma_{\gamma(0)},\Omega)\rightarrow 
 (h_{\gamma(1)}\backslash \Sigma_{\gamma(1)},\Omega)\]
is a symplectomorphism possibly not everywhere defined.
 \end{enumerate}
 \end{lemma}
 \begin{proof} This a generalization of  lemma 1.13 in \cite{Se03} for 
local symplectic forms which are of type $(1,1)$ at the origin; also -an very important for our applications- 
smooth dependence on the  symplectic form and curve $\gamma\subset \mathbb{C}$ is proved.

Points 1 and 2 are a  straightforward calculation. Point 3 is also elementary once we use Taylor expansions at the origin.

Point 3 implies that $0\in \mathbb{C}^{n+1}$ is a hyperbolic singular point for  $Y$ (see \cite{MP82} 
for basic theory on dynamical systems).
Let $W^s(Y)$ denote the stable manifold. Point 2 implies that
$[0,r_0)\subset h(W^s(Y))$ for some $r_0>0$, and that for any $r\in (0,r_0)$ the intersection $h_r\cap W^s(Y)$
is transverse. Since $\Sigma_{r}:=h_r\cap W^s(Y)$ is a hypersurface of $W^s(Y)$ transverse to $Y$, it is diffeomorphic
to a sphere. More precisely, the stable manifold theorem gives a parametrization
$\Psi^{\mathrm{st}}\colon B^{n+1}\rightarrow  W^s(Y)$ of a neighborhood of the origin inside $W^s(Y)$
 which is unique up to isotopy and the action of
$O(n+1)$, the latter associated to the choice of an orthonormal
basis of the tangent space of $W^s(Y)$ at the origin; such parametrization 
induces a parametrization $l\colon S^n\rightarrow \Sigma_r$ unique up to isotopy and the action of
$O(n+1)$.

That $\Sigma_r$ is Lagrangian follows from point 2,
exactly as in the proof  of  lemma 1.13 in \cite{Se03}.

The result for any other point $z$ and  a curve $\gamma$ joining it to the
origin follows from the previous ideas applied to the Hamiltonian
of $-\mathrm{Im}(F\circ h)$, where $F\colon \mathbb{C}\rightarrow \mathbb{C}$ is a
diffeomorphism fixing the origin which sends $\gamma$ to $[0,r]$,
for some $r\in \mathbb{R}\backslash \{0\}$.

If $\Omega_u$ is a smooth family, then the stable manifold theorem with parameters
(the proof of theorem 6.2 in \cite{MP82}, chapter 2, is seen to depend smoothly on
parameters) gives parametrizations
$\Psi^{\mathrm{st}}_u\colon B^{n+1}\rightarrow  W^s(Y_u)$
of neighborhoods of 0 inside the corresponding stable manifolds. This induces
 a smooth family of parametrizations of the Lagrangian spheres

\begin{equation*}
l_u\colon S^n\rightarrow \Sigma_{u,r}.
\end{equation*}

Clearly there is also smooth dependence  on the path
$\gamma$ if we choose diffeomorphisms $F_\gamma\colon \mathbb{C}\rightarrow\mathbb{C}$ with such
dependence.

Parallel transport is not defined for  points in $h_{\gamma(0)}$ which  converge
to the singular point $0\in \mathbb{C}^{n+1}$, which by definition are the Lagrangian sphere $\Sigma_{\gamma(0)}$. 
Parallel transport  may send points of $h_{\gamma(0)}\backslash \Sigma_{\gamma(0)}$ away from $W$. 
For those points which do not leave $W$, which at least are those close enough to $\Sigma_{\gamma(0)}$,
parallel transport is well known to be a symplectomorphism, and this finishes the proof of the lemma.
\end{proof}

\subsubsection{The shape of the symplectic $(n+1)$-handle.}\label{sssec:symp}
A parametrized sphere $L\subset M$ together with a framing determine a diffeomorphism
 $\phi\colon H\rightarrow S^n\times \overline{B^{n+1}}(1)$, where $H$
is a compact neighborhood of $L$ and   $S^n\times \overline{B^{n+1}}(1)$ is seen as a
 subset of the boundary of the standard $(n+1)$-handle
 \[\overline{B^{n+1}}(1)\times \overline{B^{n+1}}(1)\subset \mathbb{R}^{n+1}\times \mathbb{R}^{n+1}=\mathbb{C}^{n+1}.\]
 The diffeomorphism determines the manifold with corners
 \[M\#_\phi \overline{B^{n+1}}(1)\times \overline{B^{n+1}}(1).\]
A way to smoothen the corners uses the gradient flow $Y$ of $-\mathrm{Re}h$ with respect to the Euclidean metric: let us consider
a function
\[f\colon H\backslash L\cong S^n\times \overline{B^{n+1}}(1)\backslash S^n\times \{0\}\rightarrow  \mathbb{R}^+\]
 supported in the interior of $H$,
and such that near the attaching sphere $L\cong S^n\times \{0\}$ its value 
is the time needed to flow from $H$ to a neighborhood 
of  $L'\cong \{0\}\times S^n\subset  \overline{B^{n+1}}(1)\times S^n$. Then
 \[M'=M\backslash H\cup \Phi_1^{fY}(H\backslash L)\cup L'\]
is a smoothening of the new boundary of the cobordism.
Actually, one equally thinks of using as modified handle the region bounded by $H$ and $\Phi_1^{fY}(H\backslash L)\cup L'$.

\begin{figure}[h]
\centering
\includegraphics[height=5cm]{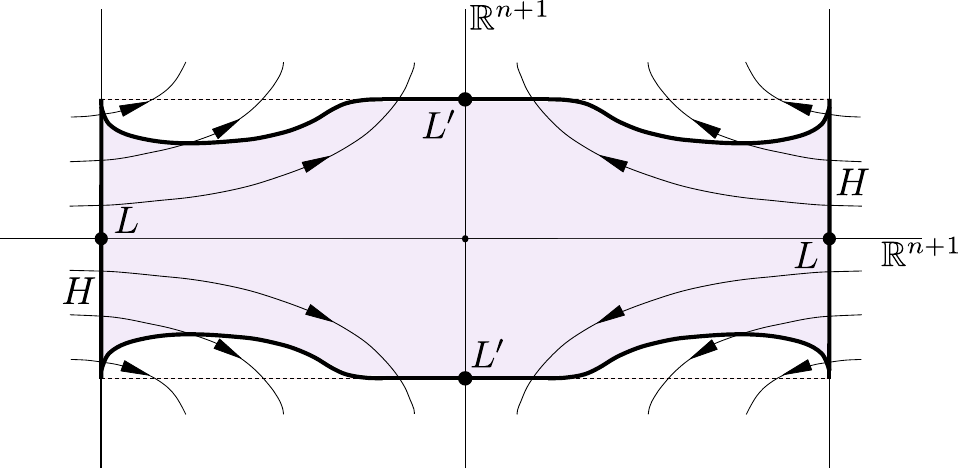}
\caption{The modified handle is the shaded region, which is everywhere transverse
 to the gradient flow lines. The doted segments are part of the boundary
of the standard handle with corners.}
\label{fig:handle}
\end{figure}

We now proceed to define smoothenings of the standard $(n+1)$-handle using $Y_2$, which is our symplectic 
replacement for the  gradient flow of $-\mathrm{Re}h$. For the sake 
of flexibility in the definition of Lagrangian surgery we make the construction depend on a small enough  parameter $r>0$.

We start by introducing some notation: the complex coordinate of $\mathbb{C}$ is $z=(x,y)$. For any $r,a,b\in \mathbb{R}$, we
let $y_r(a,b), x_r(a,b)\subset \mathbb{C}$ be the ``vertical'' and ``horizontal''
segments joining the points $(r,a)$ and $(r,b)$, and $(a,r)$ and $(b,r)$
respectively.

Let us consider $r_0>0$
small enough so that the neighborhood $(W,\Omega_2)$ of $0\in \mathbb{C}^{n+1}$ contains all Lagrangian spheres
$\Sigma_{r},\, r\in  [-r_0,0)\cup (0,r_0]$, described in point 4 in lemma \ref{lem:stb}.
We fix $\epsilon>0$ small enough and define for all $r\in (0,r_0]$ 
\[ H_{2,r}:=h^{-1}(y_r(-\epsilon,\epsilon)),\,{H}_{2,-r}:=h^{-1}(y_{-r}(-\epsilon,\epsilon)).\]
 By point 2 in lemma \ref{lem:stb}
 $Y_2\pitchfork H_{2,r}, {H}_{2,-r}$,
 so both hypersurfaces inherit 2-calibrated foliations.
 By definition of the symplectic connection, the leaves of these 2-calibrated
foliations are exactly the symplectic fibers of $h\colon H_{2,r}\rightarrow y_r(-\epsilon,\epsilon)$ and
 $h\colon {H}_{2,-r}\rightarrow y_{-r}(-\epsilon,\epsilon)$. 

The Lagrangian sphere $\Sigma_r$ is going to be the attaching sphere of the $(n+1)$-handle,
 and therefore we need to specify an isotopy 
class of parametrizations (its framing is the
Lagrangian framing): if $\Omega_2=\omega_{\mathbb{R}^{2n+2}}=\sum_{i=1}^{n+1} dx_i\wedge dy_i$,
then the Lagrangian sphere over $(r_0,0)$ is the sphere  radius $\sqrt{r_0}$ in the coordinates $x$
 \[\{(x,0)\in \mathbb{R}^{2n+2}\,|\, x_1^2+\cdots +x_{n+1}^2=r_0\}.\]  Remark that the
 Lagrangian framing is the standard framing. The subset of forms of type $(1,1)$ at the origin is convex
and hence connected (the symplectic condition holds for the segment close enough to the origin).
 We choose any path $\zeta$ connecting 
$\omega_{\mathbb{R}^{2n+2}}$ to $\Omega_2$, and lemma \ref{lem:stb} with parameter space $\zeta$ allows us to transfer the
canonical parametrization of the sphere of radius $\sqrt{r_0}$ to a parametrization  $l$ of $\Sigma_{r_0}$.
 This  completely determines the isotopy class of $l$. 

To connect the hypersurfaces $H_{2,r}$ and
${H}_{2,-r}$ we want a careful parametrization of a 
 neighborhood of $\Sigma_r$ inside $H_{2,r}$, $r\in (0,r_0]$. Let us extend 
the parametrization of  $\Sigma_{r_0}$ to a neighborhood of $\Sigma_{r_0}$  inside its leaf
\[\varphi_{r_0}\colon (U,\Omega_{\mathcal{F}})\rightarrow (T(\lambda),d\alpha_\mathrm{can}).\]
Parallel transport over the horizontal segment  $x_0(r_0,r)$ induces a parametrization of
 a neighborhood of $\Sigma_r$ inside its leaf
\begin{equation}\label{eq:sympnbhd}
\varphi_r:=\varphi_{r_0}\circ \rho_{x_0(r,r_0)}\colon ( \rho_{x_0(r,r_0)}^{-1}(U),\Omega_{\mathcal{F}})\rightarrow
 (T(\lambda),d\alpha_\mathrm{can}), \, r\in (0,r_0].
\end{equation}
We define $T_r(\lambda):=\varphi_r^{-1}(T(\lambda))$, $ r\in (0,r_0]$. 

Let $R_r$ be the (negative) Reeb vector field on  $H_{2,r}$ determined by the equality  
$h_*R_r=\tfrac{\partial}{\partial y}$, $r\in (0,r_0]$. The neighborhood
of $\Sigma_r$ inside $H_{2,r}$ that we are going to consider is  $T_r(\lambda,\epsilon)$, 
defined as in equation (\ref{eq:not})  using 
on $H_{2,r}$ the flow of $R_r$. In fact we redefine $H_{2,r}:= T_r(\lambda,\epsilon)$, $ r\in (0,r_0]$.

Let
\begin{equation}\label{eq:cutoff}
 f_r\in C^{\infty}(T_r(\lambda,\epsilon)\backslash \Sigma_r,\mathbb{R}^+)
 \end{equation}
have the following properties:
\begin{itemize}
\item The support of $f_r$ is contained in the interior of
$T_r(\lambda,\epsilon)$.
\item The time 1 flow of  $f_r Y_2$ sends  $T_r(\lambda/2,\epsilon/2)\backslash
\Sigma_r$ into
$H_{2,-r}$.
\end{itemize}
We use the  hypersurface
\begin{equation}\label{eq:hypers}
H_{2,r}^{\mu_L}:=\varPhi^{f_r Y_2}_1(T_r(\lambda,\epsilon)\backslash
\Sigma_r)\cup\Sigma_{-r}
\end{equation}
to define the handle $P_{2,r}$ as the compact domain of $\mathbb{C}^{n+1}$
bounded by $H_{2,r}^{\mu_L}$ and $H_{2,r}$. The new boundary of the cobordism is 
\[M^{\mu_L}=(M\backslash H_{2,r}) \cup H_{2,r}^{\mu_L}.\]

\subsubsection{Lagrangian surgery.}\label{ssec:lag}
\begin{proposition}\label{pro:lagattach}  Any parametrized Lagrangian sphere
 $L\subset (M^{2n+1},\mathcal{F},\omega)$, $n>1$, determines symplectic elementary cobordisms $(Z,\Omega)$ carrying a symplectic
vector field transverse to the boundary, which induce 2-calibrated foliations
 \[(M,\mathcal{F},\omega),\, (M^{\mu_L},\mathcal{F}^{\mu_L},\omega^{\mu_L}).\]
\end{proposition}
\begin{proof} Any form  of type $(1,1)$ at the origin endows the $(n+1)$-handle $P_{2,r}$ with
 a symplectic structure $\Omega_2$. The Hamiltonian vector field $Y_2$ is transverse 
to $\partial P_{2,r}$ and determines a parametrized Lagrangian sphere $\Sigma_{r}$. The parametrized
 Lagrangian sphere $\Sigma_r$ with its
Lagrangian framing is isotopic to the standard sphere with its standard framing.  Therefore applying
 corollary \ref{cor:cosympglue} produces
 the  elementary cobordism $Z$. Moreover, it gives rise to a symplectic structure $\Omega$ and a symplectic
 vector field $Y$ transverse to $\partial Z$ which
induce a 2-calibrated foliation on $\partial Z=M\coprod M^{\mu_L}$.
By construction we recover $(\mathcal{F},\omega)$ on $M$ and obtain
 $(\mathcal{F}^{\mu_L},\omega^{\mu_L})$ on $M^{\mu_L}$ which coincides
with  $(\mathcal{F},\omega)$ away from a neighborhood of $L$.
\end{proof} 

\begin{definition}\label{def:lagsurg} Let $L\subset (M^{2n+1},\mathcal{F},\omega)$, $n>1$, be a 
parametrized Lagrangian sphere. We define 
Lagrangian surgery along $L$ as any of the 2-calibrated foliations  $(M^{\mu_L},\mathcal{F}^{\mu_L},\omega^{\mu_L})$
 in proposition \ref{pro:lagattach},
 obtained as the new boundary component of the symplectic elementary cobordism which amounts to
attaching a symplectic $(n+1)$-handle as described in  \ref{sssec:shape}, \ref{sssec:symp} 
 to the trivial symplectic cobordism determined by  $(M,\mathcal{F},\omega)$.
\end{definition}

\begin{remark}\label{neglag} Instead of gluing the $(n+1)$-handle
 to the trivial cobordism we can proceed the other way around.
 This amounts to reversing the co-orientation on $(M,\mathcal{F},\omega)$ and hence
  considering on the $(n+1)$-handle the opposite symplectic vector field
$\mathrm{Im}h$.
    Actually, we can do things in an equivalent way: on the
     $(2n+2)$-dimensional $(n+1)$-handle we can use as attaching sphere
     $\Sigma_{-r}$ instead of $\Sigma_r$, $r>0$ (and also choosing
     an appropriate shape for the handle). We go from this second point
      of view to the first one  by using  the symplectic transformation
       $(z_1,\dots,z_{n+1})\mapsto (-iz_1,\dots,-iz_{n+1})$. It can be checked 
        that the new boundary is a 2-calibrated foliation
\begin{equation}\label{eq:neglag}
(M^{-\mu_L},\mathcal{F}^{-\mu_L},\omega^{-\mu_L}).
 \end{equation}
 Surgery along $L$   with framing $\mu_{L^-}$ gives (\ref{eq:neglag}) with opposite orientation.
 \end{remark}
\subsubsection{Independence on choices.}\label{sssec:indep}

In the construction of  $ (M^{\mu_L},\mathcal{F}^{\mu_L},\omega^{\mu_L})$ there are several choices 
 both in the  symplectic handle and in the trivial cobordism, which in principle may result into 
non-equivalent  2-calibrations $\omega^{\mu_L}$.
 The choices in the symplectic handle are 
the symplectic form  $\Omega_2$, the parameter $r\in (0,r_0]$ 
($r_0$ itself depends on $\Omega_2$), the function $f_r$ (this including the choice of $\epsilon>0$)  and the
 parametrization $\varphi_{r_0}$. Choices in the trivial
cobordism correspond to choices in  $H_1$. There, we have a fixed
 $l^{-1}\colon L\rightarrow S^n$ and we choose 
an extension  $\varphi\colon (U,\omega_\mathcal{F})\rightarrow (T(\lambda),d\alpha_\mathrm{can})$ and a 
Reeb vector field  $R_1$. When applying corollary \ref{cor:cosympglue} to construct the
elementary cobordism $Z$, 
the choice of extension $\varphi_{r_0}$ is absorbed into  the choice of extension $\varphi$.

In theorem   \ref{thm:equiv} we will show that for all $r>0$ small enough Lagrangian surgery produces
a 2-calibrated foliation equivalent to generalized Dehn surgery. Since according to proposition \ref{pro:gDehn} 
generalized Dehn surgery is independent
of the extension $\varphi$ and of the Reeb vector field, we just need to prove independence
of Lagrangian surgery on the function $f_r$ and the parameter $r$. Note that these two choices do
 not matter for the diffeomorphism
type of $(M^{\mu_L},\mathcal{F}^{\mu_L})$. The key technical result that provides the required flexibility 
in our  Poisson setting, is 
an extension result for  symplectomorphisms (lemma \ref{lem:inter}).

Let us first address the case when all choices are the same except for the functions 
$f_r,f'_r$ in equation (\ref{eq:cutoff}). They give rise to two hypersurfaces 
 $H_{2,r}^{\mu_L}(f_r), H_{2,r}^{\mu_L}(f'_r)$ as described in equation (\ref{eq:hypers}), 
transverse to $Y_2$ and matching near their boundary
 and near $\Sigma_{-r}$. Following the flow lines of $Y_2$ defines a compactly supported diffeomorphism
 from $H_{2,r}^{\mu_L}(f_r)$ to  $H_{2,r}^{\mu_L}(f'_r)$. The diffeomorphism 
is a Poisson equivalence because by construction it is symplectic parallel transport
over horizontal segments. Therefore the extension by the identity is a Poisson equivalence 
between the 2-calibrated foliations associated to $f_r$ and $f'_r$.
The general position argument used in the proof of theorem \ref{thm:calibnorsum} implies that this is in fact an equivalence
of 2-calibrated foliations.

The case where the only different choice is $r<r'$ is more delicate. 
We want to  construct a Poisson equivalence
\[\phi\colon  (M^{\mu_L},\mathcal{F}^{\mu_L},\omega_{r}^{\mu_L})\rightarrow 
 (M^{\mu_L},\mathcal{F}^{\mu_L},\omega_{r'}^{\mu_L})\]
which extends the identity map in the complement of $H_{2,r}^{\mu_L}\subset M^{\mu_L} $. 
Let us define $\phi_1\colon H^{\mu_L}_{2,r}\rightarrow H_{2,r'}^{\mu_L} $ to  be the map given by the flow lines of $Y_2$, 
which we just saw corresponds to symplectic parallel transport over
 horizontal segments. It is well defined near $\Sigma_{-r}$
 because for points in $\Sigma_{-r}\subset H^{\mu_L}_{2,r}$
 we make parallel transport over the segment $x_0(-r,-r')$, which does not contain the origin.

We need to introduce the following  annular subsets around  the Lagrangian sphere $\Sigma_r$, $r\in  (0,r_0]$:
\[A_r(\lambda,\lambda'):=T_r(\lambda)\backslash \mathrm{int}T_r(\lambda'), \,\lambda>\lambda'> 0,\]
\[A_r(\lambda,\lambda',\epsilon,\epsilon'):=T_r(\lambda,\epsilon)\backslash \mathrm{int}T_r(\lambda',\epsilon'),\,,
\lambda>\lambda'>0,\,\epsilon>\epsilon'>0.\]
The boundary of an annular subset is made of an inner and and outer connected component, 
according to their distance to the Lagrangian sphere.

\begin{figure}[h]
\centering
\includegraphics[height=5cm]{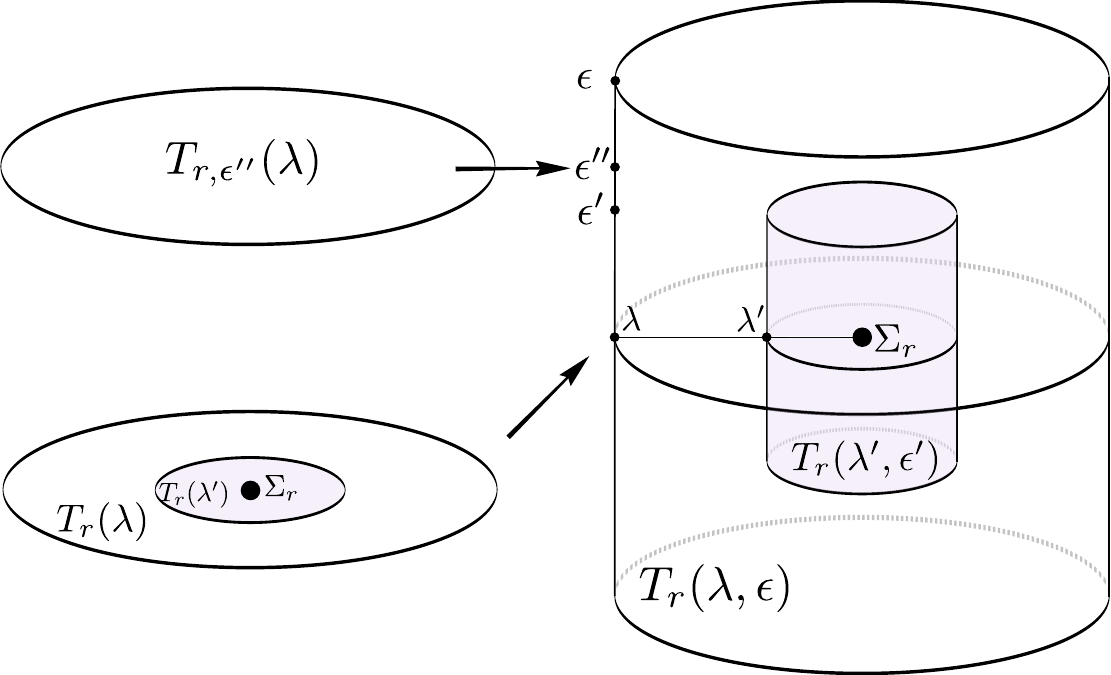}
\caption{In the r.h.s. appear the neighborhoods $T_r(\lambda,\epsilon)$ and $T_r(\lambda',\epsilon')$ (shaded)
 of the Lagrangian sphere $\Sigma_r$. Horizontal slices correspond to intersections with leaves of the foliation. In the l.h.s. 
appear the slice $t=0$, which intersects both $T_r(\lambda,\epsilon)$ and $T_r(\lambda',\epsilon')$, and the slice $t=\epsilon''$
which does not intersect $T_r(\lambda',\epsilon')$.}
\label{fig:annularshadedthin}
\end{figure}

Let $\lambda',\epsilon'>0$ be such that the support of $f_r$ (respectively $f_{r'}$) does not 
intersect  $A_r(\lambda,\lambda',\epsilon,\epsilon')$ (respectively $A_{r'}(\lambda,\lambda',\epsilon,\epsilon')$).
Therefore $A_r(\lambda,\lambda',\epsilon,\epsilon')\subset H^{\mu_L}_{2,r}\cap H_{2,r}$
 (respectively $A_{r'}(\lambda,\lambda',\epsilon,\epsilon')\subset H^{\mu_L}_{2,r'}\cap H_{2,r'}$) and on 
 $A_r(\lambda,\lambda',\epsilon,\epsilon)$ 
\begin{equation}\label{eq:interpa}
\phi_1(p)=\rho_{x_{y(h(p))}}(r,r')(p).
 \end{equation}
Note that $\phi_1$ does not extend to the identity map on 
 \[M^{\mu_L}\backslash T_r(\lambda',\epsilon')\subset M\rightarrow M^{\mu_L}\backslash T_r(\lambda',\epsilon')\subset M.\]
The problem is that according to the
 parametrizations of 
$T_r(\lambda,\epsilon)$ and $T_{r'}(\lambda,\epsilon)$ described in the paragraph following equation 
(\ref{eq:sympnbhd}), the identity map corresponds to
\begin{equation}\label{eq:interpb}
 \phi_2(p):=\rho_{y_{r'}(0,y(h(p))}\circ \rho_{x_0(r,r')}\circ \rho_{y_r(y(h(p)),0)}(p).
\end{equation}
In addition $\phi_1$ may not be everywhere defined since $\phi_1(A_r(\lambda,\lambda',\epsilon,\epsilon))$
 can fail to be contained in $A_{r'}(\lambda,\lambda',\epsilon,\epsilon')\subset H_{2,r'}\cap H_{2,r}^{\mu_L}\subset M^{\mu_L}$.

Let us assume the existence of $[\lambda_1,\lambda_1']\subset [\lambda,\lambda']$ and 
\[\phi_3\colon A_r(\lambda_1,\lambda'_1,\epsilon,\epsilon')\rightarrow A_{r'}(\lambda,\lambda',\epsilon,\epsilon')\]
a Poisson diffeomorphism onto its image, which  equals $\phi_1$ (respectively $\phi_2$) near the inner 
(respectively outer) boundary of 
$ A_r(\lambda_1,\lambda'_1,\epsilon,\epsilon')$. Then 
 \begin{equation*}
 \phi:=
 \begin{cases}
  \phi_1 & \text{in}\;\; H_{2,r}^{\mu_L}\backslash  A_r(\lambda_1,\lambda'_1,\epsilon,\epsilon'),\\
 \phi_3 & \text{in}\;\; A_r(\lambda_1,\lambda'_1,\epsilon,\epsilon'),\\
\mathrm{Id} & \text{in}\;\; M^{\mu_L}\backslash  (H_{2,r}^{\mu_L}\cap T_r(\lambda_1,\epsilon))
 \end{cases}
 \end{equation*}
is clearly an equivalence between $  (M^{\mu_L},\mathcal{F}^{\mu_L},\omega_{r}^{\mu_L})$    
  and  $ (M^{\mu_L},\mathcal{F}^{\mu_L},\omega_{r'}^{\mu_L})$.

The construction of $\phi_3$ requires the following basic result on extension
 of symplectic transformations, which is going to be also crucial to prove
the equivalence of Lagrangian and generalized Dehn surgery.
\begin{lemma}\label{lem:inter} Let 
$\varsigma_j\colon A(\lambda,\lambda')\subset (T(\lambda),d\alpha)\rightarrow (T^*S^n,d\alpha)$, $j=1,2$, $n>1$,
 be  symplectic diffeomorphisms onto their image with the following properties:
\begin{enumerate}
\item  There exists   $[\lambda_1,\lambda'_1]\subset [\lambda,\lambda']$ such that
 $\sigma_1:=\varsigma_2^{-1}\circ \varsigma_1$ is defined on $A(\lambda_1,\lambda'_1)$ and 
 there exists $\sigma_s\colon  A(\lambda_1,\lambda'_1)\rightarrow (T^*S^n,d\alpha)$, $s\in[0,1]$,
an isotopy connecting the identity to $\sigma_1$ and satisfying  
$\sigma_s(A(\lambda_1,\lambda'_1))\subset A(\lambda,\lambda')$  for all $s\in [0,1]$.
 \item They isotopy  $\sigma_s$ is Hamiltonian.
\end{enumerate}
Then there exists $\varsigma\colon  ( A(\lambda,\lambda'),d\alpha)\rightarrow (T^*S^n,d\alpha)$ 
 a symplectic diffeomorphism onto its image, which coincides
with $\varsigma_1$ (respectively $\varsigma_2$) near the inner (respectively outer) boundary 
of $A(\lambda,\lambda')$; moreover, if the $C^0$-norm of $\sigma_s$ is small enough,
then $\varsigma$ sends $A(\lambda_1,\lambda'_1)$ into $A(\lambda,\lambda)$.

In case $\varsigma_j$, $j=1,2$, the radii $\lambda,\lambda'$, the isotopy $\sigma_s$ and
 the symplectic form $d\alpha$ depend on a smooth parameter, 
 $\varsigma$ can be arranged to depend smoothly on the parameter.

\end{lemma}
\begin{proof}
Let us define
 \[V=\bigcup_{s\in [0,1]}\sigma_s(A(\lambda_1,\lambda_1'))\times \{s\}\subset T^*S^n\times [0,1].\] Its inner (respectively outer)
 boundary  is by definition the 
union of the inner (respectively outer) boundaries of $\sigma_s(A(\lambda_1,\lambda_1'))$.

Condition 1 implies $V\subset A(\lambda,\lambda')\times [0,1]$.
Let $X$ be the  vector field on $V$ whose flow after projection on $T^*S^n\times \{0\}$ gives the isotopy $\sigma_s$. 
Let $\beta_s=i_{X_s}d\alpha$. Since $\sigma_s$ is Hamiltonian
there exists a (time dependent) Hamiltonian $F\in C^\infty(V)$ such that $dF_s=\beta_s$ and $F_0=0$. 

Because $\sigma_s$ in an isotopy, for each $s\in [0,1]$ the subset  
\[A(\lambda,\lambda')\backslash \sigma_s(A(\lambda_1,\lambda_1'))\]
has an outer connected component $C_{o,s}$ (containing the outer boundary of $A(\lambda,\lambda')$)  and an inner
connected component $C_{i,s}$. We define
 \[\tilde{V}=\bigcup_{s\in [0,1]}(\sigma_s(A(\lambda_1,\lambda_1'))\cup C_{o,s})\times \{s\}\subset  A(\lambda,\lambda')\times [0,1].\]
 Let $\tilde{F}\in C^\infty(\tilde{V})$ be a function which coincides with $F$ near the inner boundary of $V$, is supported
inside $V$ and  vanishes for $s=0$. 
Then the time 1 flow of the path of Hamiltonian vector fields of  $\tilde{F}$ composed with $\varsigma_2$,
 is a symplectomorphism which coincides with  $\varsigma_1$ (respectively $\varsigma_2$) near the
 inner (respectively outer) boundary of
  $A(\lambda,\lambda'_1)$. The lemma is proved once we extend the
 symplectomorphism to $A(\lambda,\lambda')$ by using
 $\varsigma_1$ on $A(\lambda_1',\lambda')$. 

It is also clear that if the $C^0$-norm of the isotopy is arbitrarily small,
we can pick $\hat{\lambda}_1<\lambda_1$ so that $\sigma_s((A(\hat{\lambda}_1,\lambda_1'))\subset A(\lambda_1,\lambda')$, and therefore 
$\varsigma(A(\lambda_1,\lambda_1'))\subset A(\lambda,\lambda')$.
 \end{proof}
\begin{remark}\label{rem:inter} 
There is an analogous symplectic extension result when $\varsigma_j$, $j=1,2$, are defined on $T(\lambda)$. 
Under assumption 1 (with domain $T(\lambda_1)$ instead of $A(\lambda_1,\lambda_1')$), the outcome is $\varsigma$
 a symplectomorphism which matches  $\varsigma_1$ in a neighborhood of $T(\lambda')$ and $\varsigma_2$
 near the boundary of $T(\lambda)$. If the $C^0$-norm of the isotopy is small enough, then we can assume as well
$\varsigma(T(\lambda_1))\subset T(\lambda)$.
\end{remark}

We are going to apply lemma \ref{lem:inter} in several instances in which the isotopy  $\sigma_s$ 
is defined by symplectic parallel transport 
over curves $\gamma_s$. To that end
we are going to recall a straightforward result to control the $C^0$-norm of $\sigma_s$. Before that we need to
 introduce some notation. 
Given curves $\gamma_1,\dots,\gamma_n\subset \mathbb{C}$ parametrized by the interval and such that
 $\gamma_l(1)=\gamma_{l+1}(0)$, $l=1,\dots,n-1$,
  their concatenation is the piecewise
smooth curve \[ \gamma_1*\cdots *\gamma_n,\, v\in [(l-1)/n,l/n]\mapsto \gamma_l(n(v-(l-1)/n)),\,\,l=1,\dots,n-1.\]
If we speak of a  family of piecewise smooth curves, it is understood that  all the curves can
 be written as concatenation of the same  number of curves and the family is smooth
on each of the intervals.

Once we have fixed  a symplectic form $\Omega$ on  a neighborhood $W$   of the origin  which makes the fibers of the 
quadratic form $h$ symplectic,
 any piecewise smooth curve $\gamma\subset \mathbb{C}$ inside the image of $h$ 
induces by parallel transport a symplectomorphism $\rho_\gamma$, which in general is not everywhere defined on $h_{\gamma(0)}$ 
(both for points converging to the critical points and for
points escaping $W$): we just need to pull back the symplectic fibration 
$f\colon (W\backslash\{0\},\Omega)\rightarrow \mathbb{C}\backslash\{0\}$, and follow over each smooth 
piece of the curve the 1-dimensional
kernel of the closed 2-form induced on the pullback fibration. From now on and 
unless otherwise stated, by a curve  $\gamma\subset \mathbb{C}$
  we will mean a piecewise smooth curve such that
 on each smooth interval it is either constant or embedded. In this way  (i) we can define horizontal lifts 
of $\gamma$ 
without using pullback bundles, and (ii) on each smooth interval  $\gamma$ is the integral curve of a locally defined vector field.
 These two properties will make our proofs more transparent.

We also recall that $A_{r,t}(\lambda,\lambda')$, $r\in (0,r_0]$, $t\in [-\epsilon,\epsilon']$, stands for 
the time $t$ Reeb flow of $A_r(\lambda,\lambda')$, 
where the Reeb vector field is $R_r$. If we let  $\tilde{Y}$ denote the horizontal lift
of $\tfrac{\partial}{\partial y}$, then $R_r=\tilde{Y}$. Then we also define 
$A_{0,t}(\lambda,\lambda'):=\Phi_t^{\tilde{Y}}(A_0(\lambda,\lambda'))$ ($A_0(\lambda,\lambda')$ itself 
well defined because $A_{r_0}(\lambda,\lambda')\cap \Sigma_{r_0}$ is empty).

\begin{lemma}\label{lem:ode} Let $\kappa_{t,s}\subset \mathbb{C}$, $t\in [\delta,\delta']$, $s\in [0,1]$,
 be a family of loops. Let $\gamma_{t,s,l}$ be a sequence of families of loops converging to $\kappa_{t,s}$ in the $C^1$-norm 
 uniformly on $t,s$. If the horizontal lifts 
$\tilde{\kappa}_{t,s}$ starting at  $A_{r,t}(\lambda,\lambda')$ are defined for all $v\in [0,1]$ (the lift neither  converges
 to $0\in \mathbb{C}^{n+1}$ nor leaves $W$), then the following holds:
\begin{enumerate}
 \item As $l$ tends to infinity we have convergence
\[\rho_{\gamma_{t,s,l}}\overset{C^0}{\longrightarrow} \rho_{\kappa_{t,s}}\]
on $A_{r,t}(\lambda,\lambda')$ uniformly on $t,s$;
\item For any fixed $t$ if $\rho_{\kappa_{t,0}},\rho_{\gamma_{t,0,l}}$ are the identity map, $\gamma_{t,s,l}$
 does not intersect the origin and
the homotopies  $\gamma_{t,s,l}$ converge to the homotopy  $\kappa_{t,s}$ in the $C^2$-norm,
 then $\rho_{\kappa_{t,s}}$ is a Hamiltonian isotopy. 
\end{enumerate}
\end{lemma}
\begin{proof} Recall that $\rho_{\kappa_{t,s}}=\tilde{\kappa}_{t,s}(1)$. 
 Let  $K$  be the union of the horizontal lifts $\tilde{\kappa}_{t,s}$ starting at all $p\in A_{r,t}(\lambda,\lambda')$
 for all $t,s$. By
assumption $K\subset W$ is a compact subset not containing the critical point $0\in \mathbb{C}$. Then we can work inside
$U_K\subset W$ a compact neighborhood of $K$ missing the critical point, where the convergence in point 1
follows from basic ODE theory.

If $n=2$ then $\kappa_{t,s}$ may not be Hamiltonian because  $A_{r,t}(\lambda,\lambda')$ has non-trivial
first Betti number. If $\gamma_{t,s,l}$ does not contain the origin, then parallel transport cannot converge 
to the critical point $0\in \mathbb{C}^{n+1}$. It cannot
scape $W$ for connectivity reasons: for each fixed $t$ and for $l$ large enough, parallel transport  $\rho_{t,s.l}$,
 $s\in [0,1]$, is an isotopy sending  $A_{r,t}(\lambda,\lambda')$ 
inside $h_{\kappa_{t,s}(0)}\cap W$. Then it must send $T_{r,t}(\lambda)$ inside $h_{\kappa_{t,s}(0)}\cap W$.

 Because $T_{r,t}(\lambda)$ has trivial first Betti number
$\rho_{\gamma_{t,s,l}}$ is a Hamiltonian isotopy. Because convergence of the homotopies in the $C^2$-norm 
implies converges of the isotopies in the $C^1$-norm,  the closed 1-form $\beta_s$ associated to the 
isotopy $\rho_{\kappa_{t,s}}$, $s\in [0,1]$, can be $C^0$-approximated
by exact ones, and therefore it is exact and  $\rho_{\kappa_{t,s}}$ is Hamiltonian. 
\end{proof}
\begin{remark}\label{rem:ode} 
A similar convergence result holds if the horizontal lifts start at all points in $T_{r,t}(\lambda)$.
\end{remark}

We are ready to construct   $\phi_3$ on $A(\lambda_1,\lambda'_1,\epsilon,\epsilon')$ which coincides with the Poisson morphism
$\phi_1$ in (\ref{eq:interpa}) (respectively $\phi_2$  in  (\ref{eq:interpb})) near the inner (respectively outer) boundary of 
$A(\lambda_1,\lambda'_1,\epsilon,\epsilon')$.

\begin{figure}[h]
\centering
\includegraphics[height=5cm]{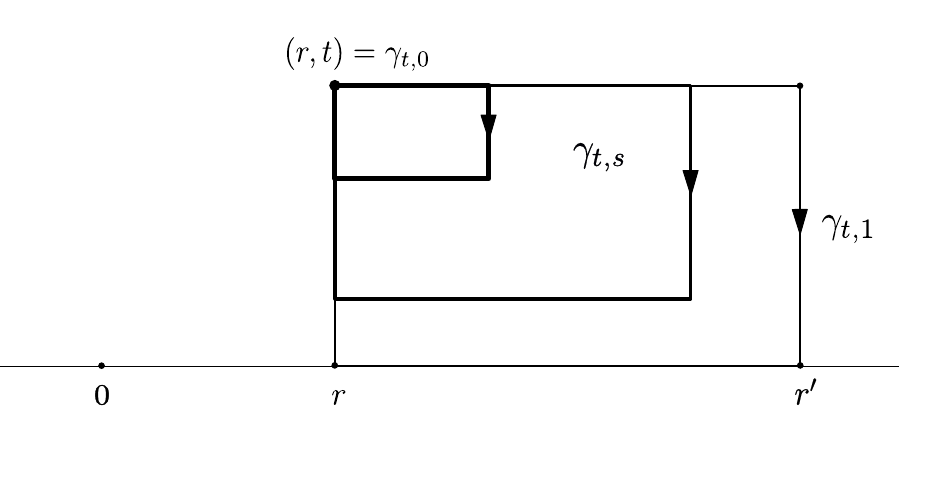}
\caption{A family of curves shrinking $\gamma_{t,1}$  the boundary of the rectangle to the vertex  $(r,t)$}
\label{fig:square}
\end{figure}

Recall that $t$ is the coordinate on the interval $[-\epsilon,\epsilon]$ and fix 
$\epsilon''\in (\epsilon',\epsilon)$. In a first stage we are going to apply lemma \ref{lem:inter} to 
the restrictions to the $t$-leaf  $\phi_{1,t},\phi_{2,t}$ with parameter
space $t\in [-\epsilon'',\epsilon'']$: 
let us define 
\[\gamma_{t,1}:=x_t(r,r')*y_{r'}(t,0)*x_0(r',r)*y_r(0,t).\]
 By equations (\ref{eq:interpa}) and (\ref{eq:interpb})
\[\sigma_t:=\phi_{2,t}^{-1}\circ \phi_{1,t}=\rho_{\gamma_{t,1}}.\]
We let $\sigma_{t,s}:=\rho_{\gamma_{t,s}}$, where
  $\gamma_{t,s}$ is a family of  curves in $\mathbb{C}$  connecting  the constant
path $(r,t)$ to $\gamma_{t,1}$, for example as depicted  in figure \ref{fig:square}.

 To get control on the $C^0$-norm of $\rho_{\gamma_{t,s}}$, 
we define $\kappa_{t,s}=y_r(t,(s-1)t)*y_r((s-1)t,t)$ and we let the family $\gamma_{t,s}$ vary with $r'$, so that
when $r'$ converges to $r$ the
 curves $\gamma_{t,s}$ converge  to the curves
 $\kappa_{t,s}$ in the $C^1$-norm.  Since $\kappa_{t,s}$ does not contain the origin and $\rho_{\kappa_{t,s}}=\mathrm{Id}$, 
by lemma 
\ref{lem:ode} if $r'$ is close enough to $r$ then $\rho_{\gamma_{t,s}}$ is as close as desired
to the identity on $A_{r,t}(\lambda',\lambda)$ in the $C^0$-norm. Remark that here we do not use the full
 power of lemma \ref{lem:ode},
as the curves $\kappa_{t,s}$ do not contain the origin.

The conclusion is that for $t\in [-\epsilon',\epsilon']$  
hypothesis 1 in lemma \ref{lem:inter}  is satisfied (it is understood that we conjugate the isotopy problem 
in $A_{r,t}(\lambda,\lambda')$ to an isotopy problem
in $A(\lambda,\lambda')$, using minus the Reeb flow for time $t$ and the chart $\varphi_r$). We can perform
 exactly the same construction
for  $|t|\in [\epsilon',\epsilon'']$ with the maps $\rho_{\kappa_{t,s}}=\mathrm{Id}, \rho_{\gamma_{t,s}}$ 
defined now on $T_{r,t}(\lambda)$, to conclude that the hypothesis 
 of remark \ref{rem:ode} is also satisfied.

Because $\kappa_{s,t}$ does not contain $0\in \mathbb{C}$, for $r'$ close enough to $r$ the isotopy
 $\gamma_{t,s}$ misses the origin,
and therefore it is a Hamiltonian isotopy. Thus we are in the hypothesis of lemma \ref{lem:inter} and 
remark \ref{rem:inter}. Inspection of the proof of lemma \ref{rem:inter} shows that the lemma and remark  can be combined
 to produce $\phi_{3,t}$, $t\in [-\epsilon'',\epsilon'']$,  depending smoothly on $t$ and 
extending  $\phi_1$ and $\phi_2$.

The extension for $|t|\in
  [\epsilon'',\epsilon]$ is straightforward: we let $\tilde{\sigma}_{t,s}$ be the isotopy corresponding
 to the Hamiltonian $\tilde{F}_t$
 in the proof of lemma \ref{lem:inter} (rather in the proof of remark \ref{rem:inter}). We have
defined  $\phi_{3,t}:=\phi_{2,t}\circ \tilde{\sigma}_{t,1}$. Let  $\beta\colon [\epsilon'',\epsilon]\rightarrow [0,1]$
be orientation reversing  and constant near the boundary. For $t\in [\epsilon'',\epsilon]$ we set $\phi_{3,t}:=\phi_{2,t}\circ \sigma_{\epsilon'',\beta(s)}$. For negative $t$
 we proceed analogously and this produces the required extension  $\phi_{3,t}$, $t\in [-\epsilon,\epsilon]$.

We have showed  that Lagrangian surgery produces equivalent 2-calibrations if $r,r'$ are close enough, 
which obviously implies the 
independence of the construction on $r\in (0,r_0]$.

\subsection{Generalized Dehn surgery is equivalent to Lagrangian surgery.}\label{ssec:poisson}
Equivalence of the two surgeries will remove all dependences appearing 
in Lagrangian surgery. The proof of the equivalence bears much resemblance to the
proof  of the independence of Lagrangian surgery on the parameter $r>0$, though it has additional technical
complications. We  give a brief overview in the following paragraphs.

To construct the equivalence between $(M^L,\mathcal{F}^L,\omega^L)$ and 
 $(M^{\mu_L},\mathcal{F}^{\mu_L},\omega^{\mu_L})$ a Poisson diffeomorphism suffices.
 The morphism is defined to be the identity away from a neighborhood of the Lagrangian spheres; 
 then -working already in the symplectic handle we used in the cobordism-  following the flow
lines of $Y_2$ extends the identity to a morphism
 \[\phi_2\colon H_{2,r}\backslash T_r(\lambda/2,\epsilon/2)\rightarrow H^{\mu_L}_{2,r}.\]
For some $\lambda'\in (\lambda/2,\lambda)$, $\epsilon'\in (\epsilon/2,\epsilon)$,
$\phi_2$ restricts to   $A_r(\lambda',\lambda/2,\epsilon',\epsilon/2)$ to parallel transport over horizontal 
segments $x_t(r,-r)$. 

Let us cut $T_{r}(\lambda',\epsilon')$ along $T_r(\lambda')$  and let
  $\chi\colon T_r^+(\lambda')\rightarrow T_r^-(\lambda')$ be conjugated to a generalized Dehn twist supported in the interior of  
$T(\lambda/2)$. We would be done if perhaps after modifying  $\phi_2$ near  the inner boundary 
of $A_r(\lambda',\lambda/2,\epsilon',\epsilon/2)$, we can extend it to a morphism 
\begin{equation}\label{eq:extdehnglue}
\phi_3\colon T_r^+(\lambda',\epsilon')\#_{\chi}T_r^-(\lambda',\epsilon')\rightarrow H_{2,-r}\subset H^{\mu_L}_{2,r}.
 \end{equation}
Equivalently we need a pair of morphisms
 $\phi^{\pm}_3\colon  T_r^{\pm}(\lambda',\epsilon')\rightarrow H_{2,-r}$ which satisfy
\begin{equation}\label{eq:gluemon}
 \phi^+_3(p)=\phi^-_3\circ \chi(p),\,p\in T_r^+(\lambda'),
\end{equation}
and which are independent of $t$ for $|t|$ small, so the induced morphism $\phi^+_3\#_\chi\phi^-_3$ is smooth.

If  $\Omega_2$  was the closed 2-form $\Omega_\tau$, then  $\chi$ can be taken
 to be 
$\rho_{\partial \bar{D}(r)}$   parallel  
transport over $\partial \bar{D}(r)$ counterclockwise. 

Consider the positive half disks  
\[\partial \overline{D}^+(r):=\{re^{i\theta\pi}\,|\, 0\leq \theta \leq 1\}\subset \mathbb{C},\]
and   set $\zeta_{t}=y_r(t,0)*\partial \overline{D}^+(r)*y_{-r}(0,t)$. Then define 
\begin{eqnarray*}
 \rho^+ \colon  T_r^+(\lambda',\epsilon') & \longrightarrow &  H_{2,-r}\\
  p &\longmapsto & \rho_{\zeta_{y(h(p))}}(p),
\end{eqnarray*}
and define $\rho^-$ on $T_r^-(\lambda',\epsilon')$ by parallel transport over the reflection of $\zeta_t$ in the $x$-axis.
Then $\rho^{\pm}$ satisfy equation (\ref{eq:gluemon}) and therefore they induce a morphism as in equation (\ref{eq:extdehnglue}).
But this morphism does not match $\phi_2$ because for the latter we do parallel transport over horizontal segments and for 
$\rho^{\pm}$ we use half disks (up to composition with vertical segments). So our problem reduces to define Poisson equivalences
on $A_r^{\pm}(\lambda',\lambda/2,\epsilon',\epsilon/2)$, which extend parallel transport over $\zeta_t$ (and its reflection 
in the $x$-axis) near the inner boundary
and parallel transport over $x_t(-r,r)$ near the outer boundary. Of course, the extensions $\phi^{\pm}_3$ have to be compatible on
$A^{\pm}_r(\lambda',\lambda/2)$ with $\chi$; because $\chi$ is supported in the interior of $T_r(\lambda/2)$ the extensions must
 coincide on $A^{\pm}_r(\lambda',\lambda/2)$.  This compatibility condition is going to follow from a careful 
 choice of the families of curves connecting  $x_r(-t,t)$ to $\zeta_t$.
Since we will be doing parallel transport near the critical point, we will need the full power
of lemma \ref{lem:ode} to argue
that we can control the norm of the isotopies we construct and hence we are in the hypothesis of the interpolation lemma. 

A further technical complication appears because the symplectic form $\Omega_2$ in the handle is different from $\Omega_\tau$. So
 the extension of parallel transport over segments and half disks has to include 
a deformation from  parallel transport with respect to  $\Omega_2$ to parallel transport with respect to $\Omega_\tau$.

\begin{theorem}\label{thm:equiv} Under the assumption $n>1$ we have equivalences of 2-calibrated foliations
\begin{equation}
\phi\colon (M^L,\mathcal{F}^L,\omega^L)\rightarrow
(M^{\mu_L},\mathcal{F}^{\mu_L},\omega^{\mu_L})
\end{equation}
for all $r>0$ small enough.
\end{theorem}

\begin{proof}
\emph{Stage 1.}
The complement $(M,\mathcal{F},\omega)\backslash T_r(\lambda,\epsilon)$ can be seen
as a subset of both $(M,\mathcal{F},\omega)$ and 
$(M^{\mu_L},\mathcal{F}^{\mu_L},\omega^{\mu_L})$.
We may assume without loss of generality that for some $\lambda'>\lambda/2,\epsilon'>\epsilon/2$, the time 1 flow of $f_rY_2$ sends
 $T_r(\lambda',\epsilon')\backslash \Sigma_r$ into  $H_{2,-r}\subset H^{\mu_L}_{2,r}\subset M^{\mu_L}$. We  define

\begin{equation*}
 \phi_0=
 \begin{cases} \mathrm{Id}  & \text{in}\;\;  M\backslash T_r(\lambda,\epsilon),\\
 \varPhi_1^{f_rY_2} & \text{in}\;\; A_r(\lambda,\lambda/2,\epsilon,\epsilon/2),
 \end{cases}
 \end{equation*}
which a Poisson morphism given on $A_r(\lambda,\lambda/2,\epsilon',\epsilon/2)$
by parallel transport over horizontal segments $x_t(-r,r)$, $t\in [-\epsilon',\epsilon']$. 

\emph{Stage 2.} 
In both $H_{2,r}$ and $H_{2,-r}$ we have Reeb vector fields $R_r$, $R_{-r}$ defined
 near $\Sigma_r$ and $\Sigma_{-r}$ respectively (they are horizontal lifts of $\tfrac{\partial}{\partial y}$). 
Their flow parametrizes the leaf spaces by $t\in [-\epsilon,\epsilon]$.
 For the purpose of checking the smoothness
of the morphism $\phi\colon M^L\rightarrow M^{\mu_L}$ in the statement of the theorem, in this stage we shall modify $\phi_0$
 near the inner boundary
 of $A_r(\lambda',\lambda/2,\epsilon',\epsilon/2)$ to make it $t$-invariant   for 
$|t|$ small (equivariant with respect to the flows of $R_r$ and $R_{-r}$).

Let $\beta\colon [-\epsilon',\epsilon']\rightarrow [-\epsilon',\epsilon']$ be an odd monotone function which is the identity
near the boundary and maps to zero exactly the interval $[-\delta,\delta]$, with $0<\delta<\epsilon/2$.
Set  
\[\zeta_{t}:=y_r(t,\beta(t))*x_{\beta(t)}(r,-r)*y_{-r}(\beta(t),t),\,\,t\in  [-\epsilon',\epsilon']\]
and define
\[\phi_1(p)=\rho_{\zeta_{y(h(p))}}(p),\,\,p\in  A_r(\lambda',\lambda/2,\epsilon',\epsilon/2),\]
which by construction is $t$-invariant for $t\in [-\delta,\delta]$.

We are going to construct 
 $\phi_2'\colon A_r(\lambda',\lambda/2,\epsilon',\epsilon/2)\rightarrow H_{2,-r}$ extending  $\phi_0$ near the outer boundary
 and $\phi_1$ near the inner boundary, 
by applying lemma \ref{lem:inter}: let  
\begin{equation}\label{eq:curveiso}\gamma_{t,s}=y_r(t,(1-s)t+s\beta(t))*x_{(1-s)t+s\beta(t)}(r,-r)*y_{-r}((1-s)t+s\beta(t),t)*x_t(-r,r),
\end{equation}
with $t\in  [-\epsilon',\epsilon']$, $s\in [0,1]$ and $r\in (0,r_0]$.  Parallel transport 
over $\gamma_{t,s}$  defined on $A_{r,t}(\lambda',\lambda/2)$ connects 
 the identity map  to $\phi^{-1}_{0,t}\circ \phi_{1,t}$. 
To estimate the $C^0$-norm of $\rho_{\gamma_{t,s}}$ we define  $\kappa_{t,s}$ by using the formula
of $\gamma_{t,s}$ in (\ref{eq:curveiso}) for $r=0$, and consider $\rho_{\kappa_{t,s}}$ with domain $A_{0,t}(\lambda',\lambda/2)$.
By construction $\rho_{\kappa_{t,s}}$ is the identity.

Let $\gamma_{t,s}'$ be the conjugation of $\gamma_{t,s}$ by $x_{t}(0,r)$ and let us consider $\rho_{\gamma'_{t,s}}$ defined on  
$A_{0,t}(\lambda',\lambda/2)$,  the same domain as for $\kappa_{t,s}$.

We construct the extension $\phi_{2,t}'$ first for the leaves in $[-\epsilon/2,\epsilon/2]$: the union of the horizontal lifts
 $\tilde{\kappa}_{t,s}$ at  $A_{0,t}(\lambda',\lambda/2)$ is exactly 
 \begin{equation}\label{eq:kompact}
K=\bigcup_{t\in [-\epsilon/2,\epsilon/2]}A_{0,t}(\lambda',\lambda/2),
  \end{equation}
a compact subset not containing the critical point $0\in \mathbb{C}^{n+1}$.
The curves $\gamma_{t,s}'$ clearly converge in the $C^1$-norm to $\kappa_{t,s}$ as $r$ goes to zero.
 Therefore by point 1 in lemma \ref{lem:ode} there is $C^0$-convergence of $\rho_{\gamma_{t,s}'}$  to the identity.

The same result holds for $\rho_{\gamma_{t,s}}$, though not automatically since parallel transport over $x_t(0,r)$ does not
send $A_{0,t}(\lambda',\lambda/2)$ diffeomorphically into $A_{r,t}(\lambda',\lambda/2)$. This is the same situation as in the proof
of independence of Lagrangian surgery on $r$. We define
\[\tau_{t,s}=x_t(0,r)*y_r((t,(s-1)t)*x_{(s-1)t}(r,0)*y_r((s-1)t,t).\]
For $r=0$ we get $x_t(0,0)*y_0(t,(s-1)t)*x_{(s-1)t}(0,0)*y_0((s-1)t,t)$. We consider parallel transport
  $\rho_{\tau_{t,s}}$ defined on $A_{r,t}(\lambda',\lambda/2)$, which for
$r=0$ is the identity. Since for $r=0$ the union of the horizontal lifts of $\tau_{t,s}$ starting at
 $A_{0,t}(\lambda',\lambda/2)$ is again $K$ in (\ref{eq:kompact}),
 by  point 1 in lemma \ref{lem:ode} we conclude that 
$\rho_{x_{t}(0,r)}(A_{0,t}(\lambda',\lambda/2))$ converges to  $A_{r,t}(\lambda',\lambda/2)$ in the
$C^0$-norm as $r$ tends to zero, and this finishes the proof of the estimate needed in point 1 of
 lemma \ref{lem:inter}   for $t\in [-\epsilon/2,\epsilon/2]$.

 For $|t|\in [\epsilon/2,\epsilon']$ the estimate holds by connectivity arguments already mentioned: the
 proof above shows that for some 
interval $[\lambda'_1,\lambda'_2]\subset [\lambda',\lambda/2]$, the isotopy  $\rho_{\gamma_{t,s}}$ sends 
$A_{r,t}(\lambda'_1,\lambda'_2)$ into $A_{r,t}(\lambda',\lambda/2)$, for $|t|\in [\epsilon/2,\epsilon']$.
 Hence it must send $T_{r,t}(\lambda_1')$ into 
$T_{r,t}(\lambda')$. 

The isotopies $\rho_{\gamma_{t,s}}$ are Hamiltonian: if $t$ is not in $[-\delta,\delta,]$, then  $\rho_{\gamma_{t,s}}$
 extends to $T_{r,t}(\lambda')$ because $\gamma_{t,s}$
 does not contain the origin. For
the remaining values of $t$ it easy to check that the homotopy  $\gamma_{t,s}$ can be approximated
 in the $C^2$-norm by a homotopy which
 does not contain $0\in\mathbb{C}$. Therefore by
point 2 in lemma \ref{lem:ode} the isotopies are Hamiltonian. Hence we can apply lemma \ref{lem:inter}
and remark \ref{rem:inter} in a compatible manner to produce
 $\phi_2'$ on $A_r(\lambda',\lambda/2,\epsilon''\epsilon/2)$, $\epsilon''\in (\epsilon/2,\epsilon')$,
 extending $\phi_0$ and $\phi_1$. For the $t$-leaves with $|t|\in [\epsilon'',\epsilon']$,
 we apply the same patching trick as in the construction of the extension
 $\phi_3$ at the end of \ref{sssec:indep}.

We define for $r>0$ small enough
\begin{equation*}
 \phi_2=
 \begin{cases} \phi_0  & \text{in}\;\;  M\backslash T_r(\lambda',\epsilon'),\\
 \phi_2' & \text{in}\;\; A_r(\lambda',\lambda/2,\epsilon',\epsilon/2), 
 \end{cases}
 \end{equation*}
which is  a Poisson morphism independent of  $t\in [-\delta,\delta]$.

\emph{Stage 3.} In this stage we cut $M$ along a neighborhood of $L$ inside its leaf $F_L$, and 
then define a Poisson morphism  which extends  $\phi_2$ in stage 2   and parallel transport
 over boundaries of half disks (``conjugated'' by vertical segments); the latter parallel transport
also includes a deformation from  $\Omega_2$ to $\Omega_\tau$.

Let us assume for the moment that $\Omega_2$ equals $\Omega_{\mathbb{R}^{2n+2}}$. 
The closed 2-forms $\Omega_\tau$ (\cite{Se03}, section 1.2) 
are written
\[\Omega_\tau= \Omega_{\mathbb{R}^{2n+2}}+d\alpha,\]
where $d\alpha$ vanishes on the tangent space to the fibers $h_z$ and is zero in a neighborhood of the union
 of stable and unstable manifold of $Y_2$
with respect to  $\Omega_{\mathbb{R}^{2n+2}}$. The first property implies that the fibers $h_z$ are symplectic. 
The second property implies that
symplectic parallel transport with respect to $\Omega_\tau$ over $x_0(r,-r)$ is defined on $T_{r}(\lambda)\backslash \Sigma_r$.

 We assume that $\alpha$
has been chosen so that  parallel transport over $\partial\bar{D}(r)$ counterclockwise is conjugated by $\varphi_r$
 to a generalized Dehn twist
supported in the interior of $T(\lambda/2)$.
 Let us define $\Omega_u=\Omega_{\mathbb{R}^{2n+2}}+ud\alpha$, $u\in [0,1]$, and let $u\colon [0,\epsilon']\rightarrow [0,1]$
be a monotone function which attains the value 0 on $[2\delta/3,\epsilon']$ and the value 1 on $[0,\delta/3]$. 

Let us consider the arcs
\[\partial \overline{D}^+_t(r):=\{(0,t)+re^{i\theta\pi}\}\subset \mathbb{C},\]
and let us  define the curves 
\[\zeta_{t}=y_r(t,\beta(t))*\partial \overline{D}^+_{\beta(t)}(r)*y_{-r}(\beta(t),t).\]
Next we cut $T_r(\lambda',\epsilon')$ along $T_r(\lambda')$ and define on  $T_r^+(\lambda',\epsilon')$
\begin{equation}\label{eq:pertwo}
\phi_3(p)=\rho_{u(y(h(p))),\zeta_{y(h(p))}}(p),
\end{equation}
which is $t$-invariant for  $t\in [0,\delta/3]$ and on $T_{r,0}(\lambda')$ is 
parallel transport over $\partial\bar{D}(r)^+$ counterclockwise with respect to $\Omega_\tau$, and therefore
 conjugated to a Dehn twist supported
in the interior of $T(\lambda/2)$. We stress that this is a Poisson morphism because the restriction
of $\Omega_u$ to fibers of $h$ is independent of $u$ (of course what changes is the symplectic connection).

We address now the construction of $\phi^+_3$ on $A_r^+(\lambda',\lambda/2,\epsilon',\epsilon/2)$ a Poisson morphisms extending
 $\phi_2$ and $\phi_3$ and $t$-invariant for  $t\in [0,\delta/3]$, using the same pattern as in stage 2.

Let us define the curves
\begin{equation}\label{eq:segdisk}
 \gamma_{t,s}=y_{r}(t,\beta(t))*x_{\beta(t)}(r,sr)*\partial \overline{D}^+_{\beta(t)}(sr)*x_{\beta(t)}(-sr,r)*y_{r}(\beta(t),t),
\end{equation}
for $t\in [0,\epsilon']$ and $s\in [0,1]$ (see figure \ref{fig:halfdisk}). We have
 \[\phi_{\gamma_{t,1}}=\phi_{2,t}^{-1}\circ \phi_{3,t},\,\,\, \phi_{\gamma_{t,0}}=\mathrm{Id}.\]

 \begin{figure}[h]
\centering
\includegraphics[height=5cm]{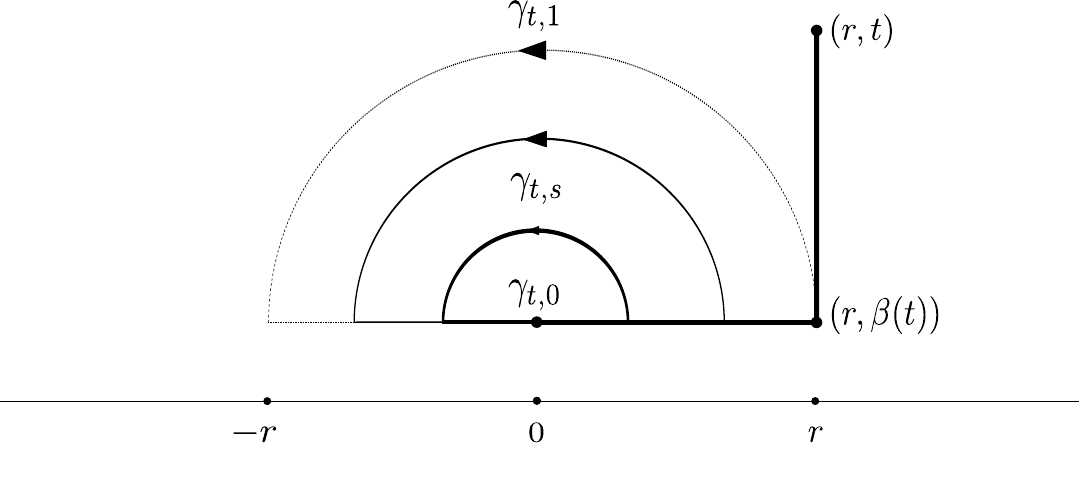}
\caption{The curves $\gamma_{t,s}$ defined in (\ref{eq:segdisk}).}
\label{fig:halfdisk}
\end{figure}

Smoothness of $\rho_{u(t),\gamma_{t,s}}$ for $s=0$ may not be evident.
\begin{lemma} The map $\rho_{u(t),\gamma_{t,s}}$ depends smoothly on $t,s$.
 \end{lemma}
 \begin{proof}
We rewrite $\rho_{u(t),\gamma_{t,s}}$ using vector fields
on $\mathbb{C}$ whose integral curves are the pieces whose concatenation defines  $\gamma_{t,s}$. 
Let
 \[X:=\frac{\partial}{\partial x}, \,Y:=\frac{\partial}{\partial y},\,
\Theta_{r,t}:=rx\frac{\partial}{\partial y}-(r(y-t))\frac{\partial}{\partial x},\, t\in \mathbb{R},\] 
be vector fields on $\mathbb{C}$.
Let $\tilde{X}_u,\tilde{Y}_u,\tilde{\Theta}_{u,r,t}\in\mathfrak{X}(W\backslash \{0\})$
be their horizontal lifts with respect to  the symplectic connection defined by $\Omega_u$.
The flows $\varPhi_l^{\tilde{X}_u}$,  $\varPhi_l^{\tilde{Y}_u}$,
$\varphi_l^{\tilde{\Theta}_{u,r,t}}$ are smooth in $u,r,t,l$. It follows that
 \begin{equation*}
\rho_{u(t),\gamma_{t,s}}=\varPhi_{t-\beta(t)}^{\tilde{Y}_{u(t)}}\circ
\varPhi_{(s+1)r}^{\tilde{X}_{u(t)}}\circ
\varPhi_\pi^{\tilde{\Theta}_{u(t),sr,\beta(t)}}\circ
\varPhi_{(1-s)r}^{-\tilde{X}_{u(t)}}
\circ\varPhi_{t-\beta(t)}^{-\tilde{Y}_{u(t)}},
\end{equation*}
and thus  $\rho_{u(t),\gamma_{t,s}}$ has smooth dependence on $t,s$. 
\end{proof}

The estimate in lemma \ref{lem:ode} is written for parallel transport
 with respect to a fixed symplectic form, but it can be checked that 
it holds true as well in case the parallel transport is with respect $\Omega_{u(t)}$.
 Let ${\kappa_{t,s}}$ be as defined in (\ref{eq:curveiso}) for $r=0$. By comparing 
 $\rho_{u(t),\gamma_{t,s}}$ with $\rho_{u(t),\kappa_{t,s}}=\mathrm{Id}$
as in the previous stage (first conjugating with $x_t(r,0)$ to have common domain, 
and then showing that the estimate holds after undoing the conjugation), 
we get control
on the $C^0$-norm of  $\rho_{u(t),\gamma_{t,s}}$ for
$r$ small enough. The isotopies  $\rho_{u(t),\gamma_{t,s}}$ are Hamiltonian since they can be 
$C^1$-approximated by Hamiltonian ones.
Given that the isotopy  $\rho_{u(t),\gamma_{t,s}}$ is $t$-invariant for $t\in  [0,\delta/3]$, choices in the 
proof of lemma \ref{lem:inter} can be done
 to obtain for all $r$ small enough an extension $\phi^+_3$ on  $A_r^+(\lambda',\lambda/2,\epsilon',\epsilon/2)$ 
 which is $t$-invariant for $t\in [0,\delta/3]$.

For $t\in [-\epsilon',0]$ we proceed as we did  for positive values,
 but using the reflection of the curves $\gamma_{t,s}$ in the $x$-axis. It is possible to 
arrange the proof of lemma \ref{lem:inter} to produce  an extension
 $\phi^-_3$ on  $A_r^+(\lambda',\lambda/2,\epsilon',\epsilon/2)$ such that
\begin{itemize}
 \item  $\phi^-_3$ is $t$-invariant for $t\in [-\delta/3,0]$;
\item  $\phi^+_3=\phi^-_3$ on $A_{r,0}(\lambda',\lambda/2)$.
\end{itemize}

Then we extend $\phi^+_3 $ to $T_r^+(\lambda',\epsilon')$ by using on  $T_r^+(\lambda/2,\epsilon/2)$ the
 same parallel transport over $\zeta_t$ as in (\ref{eq:pertwo}). 
Likewise, we extend  $\phi^-_3 $ to $T_r^-(\lambda',\epsilon')$ by using on  $T_r^-(\lambda/2,\epsilon/2)$
  parallel transport over the reflection of $\zeta_t$
in the $x$-axis.

Because   \[\phi^+_3(p)=\phi^-_3\circ \chi(p),\,p\in T_r^+(\lambda')\]
and $\phi^+_3,\phi^-_3$ are $t$-invariant for $|t|$ small, they give rise to a Poisson morphism
\[\phi_3^+\#_{\rho_{\partial\bar{D}(r)}}\phi_3^-\colon T_r^+(\lambda',\epsilon')\#_{\rho_{\partial\bar{D}(r)}}T_r^-(\lambda',\epsilon')\rightarrow H_{2,-r}.\]

The equivalence of 2-calibrated foliations for all $r>0$ small enough is

\begin{equation}\label{eq:definition2}
 \phi=
 \begin{cases}
 \phi_2 & \text{in}\;\;M^L\backslash (T_r^+(\lambda',\epsilon')\#_{\rho_{\partial\bar{D}(r)}}T_r^-(\lambda',\epsilon')),\\
 \phi_3^+\#_{\rho_{\partial\bar{D}(r)}}\phi_3^- & \text{in}\;\;  T_r^+(\lambda',\epsilon')\#_{\rho_{\partial\bar{D}(r)}}T_r^-(\lambda',\epsilon').\\
 \end{cases}
 \end{equation}

Let us now drop the assumption $\Omega_2=\Omega_\tau$.
 Let $\Omega_u$  be a path which is constant near its boundary and which connects  $\Omega_2$ to $\Omega_{\mathbb{R}^{2n+2}}$.

Recall that the neighborhoods $T_r(\lambda,\epsilon)$ have been defined with respect to $\Omega_2$.
By the parametric version of lemma \ref{lem:stb} (more specifically by the parametric
 version of the  stable manifold theorem), we have smooth parametrizations
 $\Sigma_{u,r}$, $r\in (0,r']$.  By compactness we can extend  $\varphi_{r'}$ to
parametrizations 
\[\varphi_{u,r'}\colon (T_{u,r'}(\tilde{\lambda}),\Omega_u)\rightarrow (T(\tilde{\lambda}),d\alpha_{\mathrm{can}}).\]
Then we define  the subsets $T_{u,r}(\tilde{\lambda})$   by parallel transport of
 $T_{u,r'}(\tilde{\lambda})$ over  $x_0(r',r)$ with respect to $\Omega_u$, and their associated parametrizations 
$\varphi_{u,r}:=\varphi_{u,r'}\circ \rho_{x_0(r,r')}$. The subsets $T_{u,-r}(\tilde{\lambda})\backslash \Sigma_{u,-r}$
 and their parametrizations are defined in the same manner.

We can assume without loss of generality that the inclusion  
 \begin{equation}\label{eq:incinit}
  T_{u,r}(\tilde{\lambda})\subset T_{r}(\lambda/2)
 \end{equation}
holds for all $r\in (0,r']$. This is because the parametric version of the  stable manifold theorem implies that 
\begin{eqnarray*}
 T(\lambda)\times [0,r']\times [0,1] &\longrightarrow & h_0\\
  (q,r,u)&\longmapsto & \rho_{u,x_0(r',0)}(\varphi_{u,r}^{-1}(q))
\end{eqnarray*}
is continuous, where by definition $\rho_{u,x_0(r,0)}(p)=0 \in \mathbb{C}^{n+1}$ for $p\in \Sigma_{u,r}$.

We proceed to modify both $A_r^+(\lambda',\lambda/2,\epsilon',\epsilon/2)$ and  $\phi_3$ in (\ref{eq:pertwo}) just for values of
 $t$ in $[0,\delta]$: let 
\[b_{i}\colon[0,\delta]\rightarrow [2\tilde{\lambda}/3,\lambda'],\,\, 
 b_{o}\colon [0,\delta]\rightarrow [\tilde{\lambda}/2,\lambda/2]
\]
 be  monotone increasing functions which are constant on $[0,3\delta/4]$ and near $\delta$.

Let $\upsilon\colon [0,\delta]\rightarrow [0,1]$ be a orientation reversing smooth function which
 is constant on $[0,\delta/2]$ and on $[3\delta/4,\delta]$.

We substitute  $A_{r,t}(\lambda',\lambda/2)$ by $A_{r,t}(b_o(t),b_i(t))$ and $\phi_{3,t}$ in (\ref{eq:pertwo}) by 
\[\tilde{\phi}_{3,t}:=\rho_{0,y_{-r}(0,t)}\circ \varphi^{-1}_{0,-r}\circ \varphi_{\upsilon(t),-r}\circ
\rho_{\upsilon(t),\partial \bar{D}_{\beta(t)}^+(r)}\circ \varphi^{-1}_{\upsilon(t),r}\circ \varphi_{0,r}\circ 
\rho_{0,y_r(t,0)}
\]
defined on   $A_{r,t}(b_o(t),b_i(t))$. Note that  the modification of the symplectic form only occurs
 when the domain has been modified
to $A_{r,t}(2\tilde{\lambda}/3,\tilde{\lambda}/2)$. 

By the inclusion in (\ref{eq:incinit}) the image of $A_{r,t}(2\tilde{\lambda}/3,\tilde{\lambda}/2)$ by
 $\varphi^{-1}_{\upsilon(t),r}\circ \varphi_{0,r}\circ \rho_{0,y_r(t,0)}$ is contained in  $T_{r}(\lambda/2)$. 
 If in addition $r>0$ is small
 enough, control on the $C^0$-norm of 
$\rho_{\upsilon(t),\partial \bar{D}_{\beta(t)}^+(r)}$ by $r$ implies that
\[\rho_{\upsilon(t),\partial \bar{D}_{\beta(t)}^+(r)}\circ 
\varphi^{-1}_{\upsilon(t),r}\circ \varphi_{0,r}\circ \rho_{0,y_r(t,0)}\]
sends $A_{r,t}(2\tilde{\lambda}/3,\tilde{\lambda}/2)$ into $T_{u(t),-r}(\tilde{\lambda})\backslash \Sigma_{u(t),-r}$,
 so we can compose with the chart
$ \varphi_{\upsilon(t),-r}$. Therefore $\tilde{\phi}_{3,t}$ is well defined and for $t\in [0,\delta/2]$ 
we are in the situation $\Omega_2=\Omega_{\mathbb{R}^{2n+2}}$

Then we have to choose $\Omega_{\tau}$ whose conjugation by $\varphi_{0,r}$
 ($\Omega_0=\Omega_{\mathbb{R}^{2n+2}}$) is a Dehn twist supported in the interior of $T(\tilde{\lambda}/2)$. 

We can use the same pattern to modify  the isotopy needed to apply the extension lemma with parameters (this 
time the radii of the annuli vary
with $t$). The result is an extension $\phi^+_3$ which it $t$-invariant for $t\in  [0,\delta/4]$, 
and which on $A_{r,0}(2\tilde{\lambda}/3,\tilde{\lambda}/2)$
the chart $\varphi_{0,t}$ conjugates to a Dehn twist supported on $T(\tilde{\lambda}/2)$.

As we did in the previous stage, we construct the extension $\phi^-_3$  using as domain
 and curves the reflection of the previous data in the $x$-axis.
 Then $\phi$ defined as in (\ref{eq:definition2})
 is the equivalence of 2-calibrated
foliations which proves the theorem.

\end{proof}

 \begin{remark}\label{negeq} Similarly, for $n>1$ and every $r>0$ small enough
one constructs equivalences
 \[(\phi\colon -M^{L^-},\mathcal{F}^{L^-},\omega^{L^-})\rightarrow
(M^{-\mu_L},\mathcal{F}^{-\mu_L},\omega^{-\mu_L}).\]
and
\[\phi\colon
(M^{-L},\mathcal{F}^{-L},\omega^{-L})\rightarrow(M^{\mu_{L^-}},\mathcal{F}^{
\mu_{L^-}},\omega^{\mu_{L^-}}).\]
\end{remark}

\section{Lefschetz Pencil Structures and Transverse Taut Foliations.}\label{sec:pencils}

Let $(M^{2n+1},\mathcal{F},\omega)$ be an integral 2-calibrated foliation. In this section we gather
information on the intersection of a Donaldson type submanifold with the leaves of $\mathcal{F}$ using Lefschetz
pencil structures. We also describe the relation between two Donaldson type submanifolds belonging to the same Lefschetz pencil.

We  start by saying a few words about how Donaldson type submanifolds $W$ are constructed, and how the
failure of standard Morse theoretic methods to describe the topology of $W\cap F$, $F\in \mathcal{F}$,  leads to 
the use of Lefschetz pencil structures to address this problem.

Let us fix $J$ a leafwise almost complex structure compatible with $\omega$. If $J$ is integrable then by definition
  $(M,\mathcal{F},J)$ is a Levi-flat manifold, and the line bundle $L_\omega$ whose curvature is $-2\pi i\omega$
is a positive CR line bundle. According to \cite{OS00} large powers of $L_\omega$ (suitably twisted) have plenty
of CR sections. In particular there exist CR sections leafwise transverse to the zero section
 of $L_\omega^{\otimes k}$. The zero set of any such section 
is a codimension two CR submanifold, or a divisor, intersecting $\mathcal{F}$ transversely. 

In general $J$ is not integrable. However $L_\omega^{\otimes k}\otimes \mathbb{C}^l$ has sections $s$ 
which are both close to be $J$-holomorphic in an appropriate sense and leafwise transverse to 
the zero section of $L_\omega^{\otimes k}\otimes \mathbb{C}^l$ (\cite{IM04a}, corollary 1.2). As a consequence
 $W=s^{-1}(0)$ is a 2-calibrated submanifold of 
 $(M,\mathcal{F},\omega)$ of codimension $2l$, and it is what we call a Donaldson type submanifold. The topology of $W$
and the topology of $M$ are related by a Lefschetz hyperplane type result:
the section $s$ is chosen so that   $\mathrm{log}s\bar{s}$ is  a Morse
 function. By approximate $J$-holomorphicity the index of critical points is greater than  $n-l$,
 from which the vanishing of  $\pi_i(M,W)$, $0\leq i\leq n-l-1$, follows (\cite{IM04a}, corollary 1.2).
 In particular the common zero set of $n-1$
well chosen such sections of $L_\omega^{\otimes k}$, is $W^3\hookrightarrow M$ a connected Donaldson 
type 3-dimensional submanifold.   

For any given leaf $F$, it is tempting to study the topology of $W^3\cap F$ by the same Morse theoretic methods.
It is always possible to arrange the tuple  $s=(s_1,\dots,s_{n-1})$ so that the restriction of $\mathrm{log}s\bar{s}^2$ to $F$ is 
a Morse function. The usual  Morse theoretic argument (\cite{Do96} or  proposition 2 in \cite{Au98}) implies that 
critical points have
index greater than one, and therefore if $F$ is compact  (and hence $W^3\cap F$ is compact), then $W^3\cap F$ is connected.
 If $F$ is not compact then the restriction of $\mathrm{log}s\bar{s}^2$ to $W^3\cap F$ is never proper, and it is not clear
how the information on index of critical points can be translated into topological information about $W^3\cap F$.

A second approach to study the topology of complex manifolds is via holomorphic Morse functions and Picard-Lefschetz theory.
In our setting these are Lefschetz pencil decompositions 
of $(M,\mathcal{F})$ provided  by ratios of suitable
 pairs of sections $s_1,s_2$ of $L_\omega^{\otimes k}$.  Very much as we did in the previous section with 
the complex quadratic function $h$,
 we are going to use the parallel transport associated to a Lefschetz pencil decomposition, to ``reconstruct'' a
leaf $F$ from its intersection with a regular fiber of the pencil (the previous section contains 
the analysis around a critical point of the holomorphic Morse function).  This will be enough to prove
 theorem \ref{thm:lefhyp}. Parallel transport is also the way to compare two regular fibers
 of a given Lefschetz pencil structure, showing that they differ by a sequence of generalized Dehn twists.

\subsection{Lefschetz pencil structures} We recall  the notion of Lefschetz pencil structure and the main
 existence result, and  collect some necessary results regarding the associated leafwise parallel transport.

\begin{definition} Let $x\in (M,\mathcal{F},\omega)$. A chart $\varphi_x\colon
(\mathbb{C}^n\times\mathbb{R},0)\rightarrow
(M,x)$ is compatible with $(\mathcal{F},\omega)$ if it is a foliated chart, and
$\varphi_x^*\omega$ restricted to the leaf through the origin is of type
$(1,1)$ at the origin.
\end{definition}

\begin{definition} \cite{IM04b}\label{def:pencil}
A Lefschetz pencil structure for
$(M,\mathcal{F},\omega)$ is given by a triple $(f,B,\Delta)$, where $B\subset M$ is a
codimension four 2-calibrated submanifold  and \[f\colon M\backslash
B\rightarrow \mathbb{C}\mathbb{P}^1\] is a smooth map such that:
\begin{enumerate}
\item $f$ is a leafwise submersion away from $\Delta$, a 1-dimensional manifold
transverse to
$\mathcal{F}$ where the restriction of the differential of $f$ to $\mathcal{F}$
vanishes.  The fibers of the restriction of $f$ to $M\backslash
(B\cup \Delta)$ are 2-calibrated submanifolds.
\item Around any critical point $c\in \Delta$ there exist Morse   coordinates $z_1,\dots,z_n,t$
compatible with
$(\mathcal{F},\omega)$,  and a standard complex  affine coordinate on $\mathbb{C}\mathbb{P}^1$ such that
\begin{equation}\label{eq:morsecoord}
f(z,t)=z_1^2+\cdots +z_n^2+\sigma(t),
\end{equation}
 where $\sigma\in C^\infty(\mathbb{R},\mathbb{C})$.
\item  Around any base point $b\in B$ there exist  coordinates $z_1,\dots,z_n,t$
compatible with
$(\mathcal{F},\omega)$, and a standard complex  affine coordinate on $\mathbb{C}\mathbb{P}^1$ such that
$B\equiv z_1=z_2=0$ and $f(z,t)=z_1/z_2$.
\item $f(\Delta)$ is an immersed curve in general position.
\end{enumerate}
\end{definition}

For each regular value $z\in \mathbb{C}\mathbb{P}^1\backslash f(\Delta)$, the regular fiber is the compactification
 $W_z:=f^{-1}(z)\cup B$, which is a (compact)
 2-calibrated submanifold.

\begin{theorem}[\cite{IM04b}, theorem 1.2]\label{thm:pencil} Let $(M,\mathcal{F},\omega)$ be an integral 2-calibrated
foliation  and let $e$
be an integral lift of $[\omega]$.  Then for all $k\gg 1$ there exist
Lefschetz pencils $(f_k,B_k,\Delta_k)$ such that:
\begin{enumerate}
\item The regular fibers are Poincar\'e dual to $ke$.
\item The inclusion $l_k\colon W_k\hookrightarrow M$ induces maps ${l_k}_*\colon
\pi_i(W_k)\rightarrow \pi_i(M)$ (resp. ${l_k}_*\colon
H_i(W_k;\mathbb{Z})\rightarrow H_i(M;\mathbb{Z})$) which are isomorphism for $i\leq
n-2$ and epimorphisms for $i=n-1$.
\end{enumerate}
\end{theorem}

\subsubsection{Leafwise symplectic parallel transport}\label{sssec:transport}

Let $(f,B,\Delta)$ be a Lefschetz pencil structure for $(M,\mathcal{F},\omega)$. Away from the union of 
base points and critical points $B\cup\Delta$,
the fibers of $f$ are 2-calibrated submanifolds. In particular for any point $p\notin B\cup\Delta$ this is equivalent  to
the tangent space of the leaf through $p$ and the tangent space to the fiber of $f$ through $p$ intersecting 
transversely in a symplectic
subspace.  Therefore  
 the leafwise symplectic orthogonals to the fibers
define an Ehresmann connection for $f$, which we denote by $\mathcal{H}$ and also refer to as the horizontal distribution.

The Ehresmann connection $\mathcal{H}$ is defined in the non-compact manifold $M\backslash (B\cup \Delta)$.
We are going to show that we have  good control on parallel transport near
 base points and critical points.

Let $F$ be a leaf of the foliation. Let $B_F$ (respectively $\Delta_F$) denote the codimension four 
(respectively dimension zero) submanifold
of base  (respectively critical) points in $F$. The image $f(\Delta_F)$ is a possibly countable 
collection of points in the immersed curve 
$f(\Delta)$. In particular it is easy to construct  curves $\gamma\subset \mathbb{C}\mathbb{P}^1$ which do not 
intersect $f(\Delta_F)$ (or to homotope  curves to 
avoid $f(\Delta_F)$). 

\begin{lemma}\label{lem:tamebase} Let $\gamma\subset \mathbb{C}\mathbb{P}^1$ be  curve not intersecting $f(\Delta_F)$. 
Then parallel
 transport $\rho_\gamma\colon f^{-1}_{\gamma(0)}\cap F\rightarrow  f^{-1}_{\gamma(1)}\cap F$ is 
a well defined symplectomorphism. Moreover,
it extends smoothly to a symplectomorphism 
 $\rho_\gamma\colon W_{\gamma(0)}\cap F\rightarrow  W_{\gamma(1)}\cap F$ which is the
 identity on $B_F$. In particular if $\gamma$ misses $f(\Delta)$,
it induces an  equivalence of 2-calibrated foliations
\[\rho_\gamma\colon W_{\gamma(0)}\rightarrow  W_{\gamma(1)}\] which
is the identity on $B$.
\end{lemma}
\begin{proof}
A standard procedure in this situation 
is to blow up $B$ along its leafwise almost complex normal directions.

We consider the following model for the blow up as a submanifold of $M\times \mathbb{C}\mathbb{P}^1$: 
we let $\tilde{M}$  be the union of the graph of $f$ and $B\times \mathbb{C}\mathbb{P}^1$. We need to show that  $\tilde{M}$
is a submanifold around points in  $B\times \mathbb{C}\mathbb{P}^1$. 

Around a point $b\in B$ theorem \ref{thm:pencil} provides coordinates $z_1,\dots,z_n,t$ and a standard affine coordinate
on $\mathbb{C}\mathbb{P}^1$, such that $B\cong z_1=z_2=0$ and $f=z_1/z_2$. This is equivalent to saying 
that near $b$ the graph of $f$  is given by
\begin{equation} ((z_1,\dots,z_n,t), [z_1:z_2])\subset M\times  \mathbb{C}\mathbb{P}^1.
\end{equation}
In these coordinates $\tilde{M}$ coincides with the complex blow up in the first two coordinates, and therefore it is a submanifold. 

The first projection restricts to the blow down map $\pi\colon  \tilde{M}\rightarrow M$, which is the identity away 
from $B$ and collapses
each $\{b\}\times  \mathbb{C}\mathbb{P}^1\subset \tilde{M}$ to $b\in B\subset M$. The restriction to $\tilde{M}$ of the second
projection on $M\times \mathbb{C}\mathbb{P}^1$ defines an extension of $f$
\[\tilde{f}\colon \tilde{M}\rightarrow \mathbb{C}\mathbb{P}^1.\]
Because we are blowing up directions inside leaves we have an induced foliation $\tilde{\mathcal{F}}$, and  the 
blow down map is a map of foliated
manifolds. 

The fibers of $\tilde{f}$
are transverse to $\tilde{\mathcal{F}}$, and by construction the restriction of the projection
 $\pi\colon \tilde{f}_z\rightarrow W_z$ is a diffeomorphism foliated manifolds. We
let $\tilde{F}$ denote the leaf mapping into $F$.

Let $\tilde{\omega}$ denote the pullback of $\omega$ by the blow down map. We claim that the intersection 
of the fibers of $\tilde{f}$
 with $\tilde{\mathcal{F}}$ are symplectic manifolds with respect to $\tilde{w}_{\tilde{\mathcal{F}}}$,
and therefore there is an associated leafwise Ehresmann connection which extends $\mathcal{H}$. At a point 
$p=(b, [z_1:z_2])$, say $z_2\neq 0$,
 the tangent space  $T_{[z_1:z_2]}\mathbb{C}\mathbb{P}^1\subset T_{(b,[z_1:z_2])}\tilde{F}$
 is in the kernel of $\tilde{\omega}_{\mathcal{F}}$ because the blow down map collapses the  $\mathbb{C}\mathbb{P}^1$ 
factor into the point $b$. The subspace 
$T_{(b,[z_1:z_2])}\tilde{f}\cap  T_{(b,[z_1:z_2])}\tilde{F}$ is complementary to $T_{[z_1:z_2]}\mathbb{C}\mathbb{P}^1$ and 
it is mapped isomorphically into $T_bf_{z_1/z_2}\cap T_bF$,
and the latter is symplectic with respect to $\omega_\mathcal{F}$ (alternatively, in  local coordinates about the base point
 $T_bf_{z_1/z_2}\cap T_bF$  is a complex hyperplane of $T_0\mathbb{C}^{n}$, and therefore it is symplectic
 with respect to $\omega_\mathcal{F}$  because the symplectic form has  type $(1,1)$ at the origin).
Thus, the blow down map identifies $\tilde{f}_z$ and $W_z$ as 2-calibrated foliations. 

Once we have described the kernel of $\tilde{w}_{\tilde{\mathcal{F}}}$, it is easy to see that
 the horizontal lift of $\gamma\subset  \mathbb{C}\mathbb{P}^1$ starting at $(b,\gamma(0))$ 
is exactly $(b,\gamma)$.

It is clear that parallel transport defines a Poisson equivalence. It is obviously 
an equivalence of 2-calibrated foliations because (the induced) co-orientations are preserved, and 
the 2-calibrations are restriction of the same closed 2-form on $M$.
\end{proof}

Let  $c\in \Delta$ be a critical point and let us apply  theorem \ref{thm:pencil} to construct Morse coordinates for $f$
 centered at $c$. By restricting
Morse coordinates to the leaf $F$ containing $c$, we obtain Morse coordinates for the restriction of $f$ to $F$.
 Since the restriction 
of $\omega_\mathcal{F}$ to $F$ is mapped to a symplectic form of type $(1,1)$ at the origin, leafwise parallel transport near
$c$ corresponds to parallel
transport in $\mathbb{C}^{n}\backslash \{0\}$ near $0$ for the function $h$ with respect to a symplectic form of 
type $(1,1)$ at the origin.

Let us consider the following system of neighborhoods of the critical point $0\in \mathbb{C}^n$ (\cite{Se03}, section 1.2):
 we fix the standard symplectic form $\Omega_{\mathbb{R}^{2n}}$
and define $\Sigma_z$, $z\in \mathbb{C}$, to be the Lagrangian sphere of points in $h_z$ whose parallel
 transport over the radial segment converges to the origin.
For some $r_0>0$ we fix the parametrization 
\[\varphi_{r_0}\colon (T_{r_0}(\lambda), \Omega_{\mathbb{R}^{2n}})\rightarrow (T(\lambda),d\alpha_{\mathrm{can}}).\]
For any $z\in \mathbb{C}$ small enough  we define $T_z(\lambda)\backslash \Sigma_z$ by radial parallel transport to 
the origin and then to $r_0$. Of course,
$T_z(\lambda)$ denotes the union of  $T_z(\lambda)\backslash \Sigma_z$ and $\Sigma_z$.

Then 
\begin{equation}\label{eq:parneigh}
\mathcal{T}(\lambda,r)=\bigcup_{z\in \bar{D}(r)}T_z(\lambda), \,\lambda,r>0
\end{equation}
is a system of neighborhoods of the origin. We also have the corresponding annular subsets $\mathcal{A}(\lambda,\lambda',r,r')$.

\begin{lemma}\label{lem:crittrans} Let $\Omega_u$, $u\in K$, be a compact family of symplectic forms  defined on 
a neighborhood $W$ of $0\in \mathbb{C}^n$ which make the fibers $h_z$ symplectic submanifolds. Let us fix 
 any $\lambda,r>0$,  $\lambda'\in (0,\lambda)$ and  $r'\in (0,r)$.
Then  there exists $\delta>0$ such that for any  curve $\gamma\subset \bar{D}(r')\backslash \{0\}$ having the $C^1$-norm
of $\gamma-\gamma(0)$  bounded by $\delta$,
 the horizontal lift $\tilde{\gamma}_u$ starting at any $p\in \mathcal{T}(\lambda',r')$ is contained in 
$\mathcal{T}(\lambda,r)$, for
all $u\in K$.
 \end{lemma}
 \begin{proof}
Let $\mathcal{C}$ denote the topological space of (piecewise embedded or constant) curves
contained in $\bar{D}(r)$ relative to the $C^1$-topology. 
Let us consider the subset
\[E=\{(\gamma,p,u,v)\subset\mathcal{C}\times \mathcal{A}(\lambda',\tilde{\lambda},r,r')\times K\times [0,1]\,|\, \gamma(0)=h(p)\},\]
and let us define the  continuous map
\begin{eqnarray*}
   G\colon E &\longrightarrow & W\\
  (\gamma,p,u,v) &\longmapsto & \tilde{\gamma}_u(v),
  \end{eqnarray*}
by sending a tuple to the evaluation for time $v$ of the horizontal lift of $\gamma$
with respect to $\Omega_u$ starting at $p$. The map is not everywhere defined
 since horizontal lifts may leave $W$ or converge to
the critical point, which is exactly what we want to control. However,
inside $E$ we have the subset $\mathcal{A}(\lambda',\tilde{\lambda},r,r')\times K\times [0,1]$ corresponding to constant curves.
The restriction of  $G$ to this subset is the first projection. By continuity 
an open neighborhood of $\mathcal{A}(\lambda',\tilde{\lambda},r,r')$ inside $E$ is sent into
$\mathcal{T}(\lambda,r)$. Because on $E$ we have the topology induced by the product topology, 
we conclude the existence of $\delta$ such that curves
with $||\gamma-\gamma(0)||_{C^1}< \delta$, $\gamma(0)\in \bar{D}(r')$,   have horizontal lift starting at points in
  $A_{\gamma(0)}(\lambda',\tilde{\lambda})$ contained in $\mathcal{T}(\lambda,r)$.
If in addition such a small curve does not contain  $0\in \mathbb{C}$, the connectivity argument already used a couple of times 
implies that  horizontal lifts starting at points in 
$T_{\gamma(0)}(\lambda')$ remain inside
$\mathcal{T}(\lambda,r)$, and this proves the lemma.
\end{proof}
\begin{remark}\label{rem:plaquemorse} Let $\sigma\in \mathbb{C}$ and consider the constant perturbation
of the complex quadratic form $h+\sigma$. Note that in the definition of 
$\mathcal{T}(\lambda,r)\subset \mathbb{C}^n$ for $h$ given in 
(\ref{eq:parneigh}), if we replace
radial segments joining a point $z$ to the origin by segments joining $z$ to $\sigma$, 
we  get exactly the same subset $\mathcal{T}(\lambda,r)$. Now assume that the parameter $u\in K$ in lemma
 \ref{lem:crittrans} describes 
not just the variation of symplectic forms, but a perturbation of $h$  by a constant $\sigma(u)$. Then 
lemma \ref{lem:crittrans} holds
replacing in the statement $h$ by $h+\sigma$ and the disk of radius $r'$ by the disk of radius $r'$ centered at $\sigma(u)$.  
\end{remark}

\subsection{Connected components of ${W_z}\cap F$ and leafwise parallel transport.}

Let us fix $z_0\in \mathbb{C}\mathbb{P}^1$ a regular value for $f$. Let $\gamma\subset \mathbb{C}\mathbb{P}^1$
 be a  loop  based at $z_0$
 with empty intersection with $f(\Delta_F)$. 
Lemma \ref{lem:tamebase} implies that parallel transport over $\gamma$ defines a diffeomorphism on 
$W_{z_0}\cap F$ (actually a symplectomorphism). Therefore the loop acts on connected components of 
$W_{z_0}\cap F$ and the action descends
 to $\pi_1(\mathbb{C}\mathbb{P}^1\backslash f(\Delta_F),z_0)$.
\begin{proposition}\label{pro:trivact} The action of $\pi_1(\mathbb{C}\mathbb{P}^1\backslash f(\Delta_F),z_0)$ 
on connected components of $W_{z_0}\cap F$
 is trivial.
\end{proposition}
\begin{proof}

 Let $\gamma$ be a loop based at $z_0$ and not intersecting $f(\Delta_F)$. Consider $H_s$
 a homotopy connecting $H_0=\gamma$ with the constant path $z_0$. 
We can assume without loss of generality that  $H$ misses a point of $\mathbb{C}\mathbb{P}^1$, and therefore
compose with an affine coordinate chart and work in $\mathbb{C}$.

We can assume as well that the curves $\gamma_s$ in the homotopy  coincide in the complement of an interval
$[a,b]\subset [0,1]$, and  the $C^1$-norm of $\gamma_{\mid [a,b]}-\gamma(a)$ (rescaled to have domain $[0,1]$)
  is bounded by any given $\delta>0$:
  by breaking the domain of $H$ into  $n^2$ squares of side $1/n$,  we can write $H$ as composition of $n^2$
 homotopies with the above property. It 
is possible that the starting curve $\gamma_0$ of each of the $n^2$ homotopies does intersect   $f(\Delta_F)$,
 but intersections can be removed
after a perturbation with does not affect the behavior we demand on the curves $\gamma_s$. 

We are going to control how the lifts of the curves in the homotopy behave near $\Delta_F$ using Morse coordinates, 
and away from  $\Delta_F$ using a compactness argument.
 
Let $c\in \Delta$ and let us construct Morse coordinates $z_1,\dots,z_n,t$ as in definition \ref{def:pencil}. We say
that the restriction of the coordinates to each plaque in their domain are 
 Morse coordinates for the restriction of $f$ to the plaque. In
Morse coordinates for a given plaque  the
restriction of $f$ transforms into $h+\sigma(t)$
and $\omega_\mathcal{F}$ transforms into a symplectic form of making the fibers  ${(h+\sigma(t))}_z$  
symplectic manifolds (this is because
we can construct Morse coordinates centered at any point in $\Delta$, and in  Morse coordinates on the plaque containing $c$
the perturbation $\sigma(t)$ can be taken to be trivial and $\omega_\mathcal{F}$ becomes a symplectic 
form of type $(1,1)$ at the origin).

 Let us cover $\Delta$ with a finite number of  Morse coordinates 
and let us consider their associated 1-parameter families of Morse coordinates on their plaques.
 Let us take $\lambda,r>0$ such that
$\mathcal{T}(\lambda,r)$ as defined in (\ref{eq:parneigh}) is contained in the image of Morse coordinates for each of the plaques. 
 Let us also  pick $\lambda'\in (0,\lambda)$, 
$r'\in (0,r)$ and denote by $U$ the points  in  $\tilde{M}$ whose image by at least one of the sets of  Morse coordinates
on its  plaque is contained in $\mathcal{T}(\lambda',r')$. Note that $U$ is a neighborhood of $\Delta$.

For each of our Morse coordinates its 1-parameter family of Morse coordinates on plaques is in the hypothesis of lemma
 \ref{lem:crittrans}, or rather remark \ref{rem:plaquemorse} (we assume that the parameter
space is a compact interval, and that these compact intervals cover $\Delta$).
Let $\delta_1$ be a  $C^1$-bound provided by remark \ref{rem:plaquemorse} and valid for the finite
number of 1-parameter families.

Let $V\subset V'$ be open neighborhoods of $\Delta$ in $\tilde{M}$ such that  
$V\subset \bar{V}\subset V'\subset U$.
Because $\tilde{M}\backslash V'$ is compact, there exists $\delta_2>0$ such that for any  $p\in \tilde{M}\backslash V'$ and any
 curve $\gamma\subset \mathbb{C}\mathbb{P}^1$ starting at $\tilde{f}(p)$ and
 such that  the $C^1$-norm of $\gamma-\gamma(0)$ is bounded by $\delta_2$, the horizontal
lift $\tilde{\gamma}$ 
starting at $p$ is contained in $\tilde{M}\backslash V$.

Let $\delta$ be the minimum of $\delta_1$ and $\delta_2$ as let us assume that for each $\gamma_s$ in our homotopy $H$ the
$C^1$-norm of ${\gamma_s}_{\mid [a,b]}-\gamma_s(a)$  is smaller than $\delta$. Let $\gamma_0$ and $\gamma_1$ be
the starting and ending curve of the homotopy and let $\tilde{\gamma}_0$ and $\tilde{\gamma}_1$
 be their respective horizontal lifts starting at
 $p\in W_{z_0}\cap F$. We claim that $\tilde{\gamma}_0(1)$ and  $\tilde{\gamma}_1(1)$ can be connected by a path in 
$W_{z_0}\cap F$, what
would prove the proposition.

Recall that $\gamma_s$, $s\in [0,1]$, is independent of $s$ in  the complement of $[a,b]\subset [0,1]$. 

Let us suppose that 
 $\tilde{\gamma}_0(a)\in \tilde{M}\backslash V'$. Because of the $C^1$-bound on ${\gamma_s}_{\mid [a,b]}-\gamma_s(a)$, $s\in [0,1]$,
 the horizontal lifts of  ${\gamma_s}_{\mid [a,b]}$ starting at  $\tilde{\gamma}_0(a)$ are 
defined for all $s\in [a,b]$ and belong to $\tilde{M}\backslash V$.
In particular $\tilde{\gamma}_s(b)$ is a curve in the fiber $\tilde{f}_{\gamma_0(b)}$. Since the curves ${\gamma_s}_{\mid [b,1]}$
 are all equal and avoid
$f(\Delta_F)$, we can construct the horizontal lift starting at all points in the path $\tilde{\gamma}_s(b)$. 
 What we just proved 
is that the homotopy $H_s$ has a well defined lift starting at $p$, and therefore 
$\tilde{\gamma}_s(1)$ connects  $\tilde{\gamma}_0(1)$ to $\tilde{\gamma}_1(1)$.

If $\tilde{\gamma}_0(a)\in V$ then it also belongs to $U$. If we compose with one of the fixed
Morse coordinates on the plaque $u_0$ containing $\tilde{\gamma}_0(a)$, the point  $\tilde{\gamma}_0(a)$
is sent to $q\in \mathcal{T}(\lambda',r')$. The curves ${\gamma_0}_{\mid [a,b]}$ and ${\gamma_1}_{\mid [a,b]}$ 
meet the hypothesis of
lemma  \ref{lem:crittrans} (remark \ref{rem:plaquemorse}), and therefore their horizontal lifts starting at $q$ are contained in 
$\mathcal{T}(\lambda,r)$. In particular the images $q_0$ and $q_1$  of $\tilde{\gamma}_0(b)$ and $\tilde{\gamma}_1(b)$ respectively
 belong to ${(h+\sigma(u_0))}_{\gamma_0(b)}$. All regular fibers of $h+\sigma(u_0)$ in $\mathcal{T}(\lambda,r)$ 
are diffeomorphic to $T(\lambda)$ and therefore they are connected.
Let $\zeta$ be a path in ${(h+\sigma(u_0))}_{\gamma_0(b)}$ connecting $q_0$ to  $q_1$.
Let us also
denote by $\zeta$  its image in the plaque $u_0$ by the Morse chart, which belongs to $\tilde{f}_{\gamma_0(b)}$.
Then the ending points of the lifts of ${\gamma_0}_{\mid [b,1]}$ starting at $\zeta(v)$, $v\in [0,1]$, 
connect  $\tilde{\gamma}_0(1)$ to  $\tilde{\gamma}_1(1)$.
\end{proof}

Proposition \ref{pro:trivact} is the key result to ``spread'' a connected component of $W_{z_0}\cap F$
onto $F$. Before, we need to show that $W_{z_0}\cap F$ is always non-empty. For that it suffices
to prove that $\tilde{f}(F)$  contains some regular value $z$ of $\tilde{f}$, because in that case
 we can use parallel transport over a curve joining $z$ to $z_0$ and avoiding singular values of
 $f(\Delta_F)$, to 
find points in $W_{z_0}\cap F$:  because $\tilde{M}$ is compact the regular values of $\tilde{f}$ (which are
the regular values of $f$) are an open dense subset. The subset $\tilde{f}(F)\subset \mathbb{C}\mathbb{P}^1$ has not empty interior
and therefore it contains regular values.

\begin{theorem}\label{thm:connect}
Let $(M^{2n+1},\mathcal{F},\omega)$, $n>1$, be a 2-calibrated foliation and
let $(f,B,\Delta)$ be a  Lefschetz pencil structure as in definition
\ref{def:pencil}. Then any regular fiber $W$ of the pencil intersects
every leaf of $\mathcal{F}$ in a unique connected component.
\end{theorem}
\begin{proof} Let $z_0$ be a regular value and let $F$ be a leaf. We let $C$
be a non-empty connected component of $W_{z_0}\cap F$ (it always exists since $W_{z_0}\cap F$ is non-empty).
Let us define $\Gamma_C$ to be the set of horizontal curves starting at $C$ and 
whose projection $\tilde{f}\circ \zeta$ is either an embedded curve or constant. 
We define
\[{F}_C:=\{p\in F\backslash \Delta_F\;|\; \exists\; \zeta \in
\Gamma_C, \,  \zeta(1)=p\}.\]
By construction ${F}_C$ is non-empty,
connected and contains $C$. We want to show that it is open. 

Let $p\in {F}_C$  such that the horizontal
curve $\zeta$ connects  $x\in C$ with $p$. Let us suppose that 
 the curve $\tilde{f}\circ \zeta$ is embedded (it is not constant).
 Then we can find a 1-parameter family of embedded curves $\gamma_s$, $s\in (-\epsilon,\epsilon)$, defined for time
 $v\in [0,1+\epsilon]$, and
such that the restriction of   $\gamma_0$ to $[0,1]$ is $\tilde{f}\circ \zeta$. Because $\zeta $ is contained in
 $\tilde{M}\backslash \Delta $, a compactness argument 
 implies that there exists $A$ an open neighborhood of $x$ inside $C$ and $\epsilon'>0$, such that the horizontal
 lift of ${\gamma_s}_{\mid [0,1+\epsilon']}$ starting
 at any point in $A$ exists for all $s\in (-\epsilon',\epsilon')$. It is clear that for $\epsilon'$ small enough 
\[U_p=\{y\in F\,|\, y=\tilde{\gamma}_s(v),\,\tilde{\gamma}_s(0)\in A,\, v\in (1-\epsilon',1+\epsilon'),\, s\in 
 (-\epsilon',\epsilon')\}\]
is a neighborhood of $p$ in $F_C$.

 If $\zeta$ is constant we make the previous construction for a family of radial curves starting at $z_0$, and the open
neighborhood is obtained considering horizontal lifts for time $v\in [0,\epsilon')$ starting at a neighborhood
$A$ of $p$ inside $C$ (we would be ``spreading'' the open subset $A$).  

We claim that $F_C$ does not contain a connected component of $W_{z_0}\cap F$ different from $C$. Suppose the contrary. Then 
we would have a loop $\gamma$ with a horizontal lift connecting two different connected 
components of $W_{z_0}\cap B$. Since after a small perturbation
 we can assume without loss of generality
that $\gamma$ does not intersect $\Delta_F$, this would contradict proposition \ref{pro:trivact}.

Because it is clear that any point in $F\backslash \Delta_F$ can
 be connected to $W_{z_0}\cap F$ by a horizontal curve 
lifting an embedded curve, we
conclude that connected components of  $W_{z_0}\cap F$ are in bijection with connected components of $F$, and this proves
that $W_{z_0}\cap F$ is connected. 
\end{proof}

\begin{proof}[Proof of theorem \ref{thm:lefhyp}]
  Let $(M,\mathcal{F},\omega)$ be a
2-calibrated foliation. If it
is not integral, compactness of $M$ implies that we can slightly
modify $\omega$ into $\omega'$ so that a suitable multiple $k\omega'$ defines an
integral homology class. Theorem \ref{thm:pencil} implies the
existence of a Lefschetz pencil $(f,B,\Delta)$.

Therefore by theorem \ref{thm:connect} any regular fiber
$(W,\mathcal{F}_W,k\omega'_W)$ intersects every leaf in a connected component. If the dimension of
 $W$ is  bigger than 3, we apply the same construction to $(W,\mathcal{F}_W,k\omega'_W)$. By induction we end up
with a 3-dimensional manifold with a taut foliation
$(W^3,\mathcal{F}_W)\hookrightarrow (M,\mathcal{F},\omega)$, whose intersection with every leaf of $\mathcal{F}$ is connected.
\end{proof}

\begin{proof}[Proof of theorem \ref{thm:leafhomeo}]
 Let  $l\colon W \hookrightarrow M$ be a submanifold as in theorem \ref{thm:lefhyp}. Because 
for all $F\in \mathcal{F}$ the intersection  $W\cap F$ is connected,  the map $l$ descends to a bijection of leaf spaces
\[\tilde{l}\colon
W/\mathcal{F}_W\rightarrow M/\mathcal{F}.\]

Open sets of $W/\mathcal{F}_W$ (respectively $M/\mathcal{F}$) are in one to one
correspondence with saturated open sets of $W$ (respectively $M$).

Let $V$ be an    saturated  open set of $(M,\mathcal{F})$. By definition $W\cap
V$ is an open set of $W$ which is is clearly saturated (even without
the assumption of $\tilde{l}$ being a bijection) and this shows that  $\tilde{l}$ is continuous.

Now let $V$ be an open saturated set of $(W,\mathcal{F}_W)$. We want to show
that its saturation in $(M,\mathcal{F})$, denoted by $\overline{V}^{\mathcal{F}}$, is
open to conclude that $\tilde{l}$ is open.

If $V$ is a saturated set and $x\in V$,
then $x$ is an interior point if and only if for some $T_x$ a local
manifold through $x$ transverse to the foliation, $x$ is an
interior point of $T_x\cap V$. Hence, every $x\in V$ is an interior point of
$\overline{V}^{\mathcal{F}}$. By
using the holonomy, if a point in a leaf is interior the whole leaf
is made of interior points. Since every leaf of $\overline{V}^{\mathcal{F}}$
intersects $V$, $\overline{V}^{\mathcal{F}}$ is open, and this proves the theorem.
\end{proof}

\subsection{Regular fibers and Lagrangian surgery}\label{ssec:surgfiber}

Let $W$ be a regular fiber of  a Lefschetz pencil structure for $(M,\mathcal{F},\omega)$.
Theorems \ref{thm:connect} and \ref{thm:pencil} describe the topology of  $W/\mathcal{F}_W$ and  part of
the homology and homotopy of $W$, 
in terms of the corresponding
data for $(M,\mathcal{F})$. We want to understand how  different regular fibers
 of the pencil are related as 2-calibrated foliations.

Let $z$ and $z'$ be regular values of the pencil belonging to the same connected
component of $\mathbb{C}\mathbb{P}^1\backslash f(\Delta)$, and let $\gamma$ 
be a curve in that connected component connecting $z$ to $z'$. Then
lemma  \ref{lem:tamebase} implies that
$\rho_\gamma\colon W_z\rightarrow W_{z'}$ is an
equivalence of 2-calibrated foliations.

We notice that any two arbitrary regular values $z$ and  $z'$ can
always  be joined by a curve $\gamma$ transverse to  $f(\Delta)$.

\begin{theorem}\label{thm:handle} Let $z, z'\in \mathbb{C}\mathbb{P}^1$ be two
regular values. Let
$\gamma$ be an embedded curve joining $z$ and $z'$ and transverse to
$f(\Delta)$. Then $f^{-1}(\gamma)$ is a cobordism between $W_z$ and $W_{z'}$ 
which amounts to adding one $n$-handle for each point $x\in \Delta$
such that $f(x)\subset \gamma$. More precisely, if $n>2$ and there
is only one critical point  $c\in f^{-1}(\gamma)$, then  there exists
$L\subset W_z\backslash B$ a framed Lagrangian sphere such that
$W_{z'}$ is the result of performing generalized Dehn surgery on
$W_z$ along $L$. The framed sphere are the points in $W_z$ which under parallel transport over $\gamma$ converge to $c$.
\end{theorem}
\begin{proof} 

Let  $w\in \gamma$ and $c\in \Delta$ with $f(c)=w$. Let us take Morse coordinates around $c$ and
 an affine chart on $\mathbb{C}\mathbb{P}^1$.
Let us assume for simplicity that the curve $\gamma$ in the affine chart coincides with 
a segment of the real axis. For $r>0$ small enough, we
want to construct a Poisson equivalence  
\[\phi\colon W_r\rightarrow W_{-r}\]

To that end consider the cobordism $Z=\tilde{f}^{-1}(x_0(-r,r))$, which is a manifold with boundary
 because $\tilde{f}$ is transverse to $\gamma$ 
($\mathrm{Im}\sigma'(0)\neq 0$). The attaching of the handle in this elementary cobordism 
occurs in a neighborhood of $c$, or equivalently in a neighborhood of 0 of the Morse chart, which is where we work from now on.

 We are going to arrange the current setting so that it becomes analogous to the one in theorem \ref{thm:equiv}. 

The pullback of $f$ to the
 $t$-leaf of $\mathbb{C}^n\times \mathbb{R}$ is $h+\sigma(t)$. After  reparametrization of the coordinate $t$,
we may assume without loss of generality that $\sigma(t)=(a(t),t)$.

The tangent space of $Z$ at $0\in \mathbb{C}^n\times \mathbb{R}$ is the hyperplane $t=0$.
Therefore the projection $Z\rightarrow \mathbb{C}^n$  is a local diffeomorphism with
image an open neighborhood $V$ of $0\in\mathbb{C}^n$.

We  define $\phi$ away from a neighborhood $V'\subset V$  of $0\in \mathbb{C}^n$ as follows:
\[\phi:=\rho_{x_0(-r,r)}\colon W_r\rightarrow W_{-r}.\]
We claim that it is possible to extend $\phi$ to an equivalence of 2-calibrated foliations repeating the proof of theorem
\ref{thm:equiv} with two minor modifications.

Let us define $\sigma_r(t):=(r,0)+\sigma(t)$, $r\neq 0$, $t\in [-\epsilon,\epsilon]$. Hence
the image of $W_r$ (respectively $W_{-r}$) on $V'$ is exactly
$h^{-1}(\sigma_r(t))$ (respectively $h^{-1}(\sigma_{-r}(t))$). 
Recall that Morse coordinates on the $t$-plaque send $\omega_\mathcal{F}$ to a symplectic form $\Omega_t$ which makes the fibers of
$h+\sigma(t)$ symplectic. Then it follows that the morphism
 $\rho_{x_0(-r,r)}\colon W_r\rightarrow W_{-}r$ at a point $p \in V'\cap W_r$ in the $t$-leaf 
corresponds  to  parallel transport $\rho_{t,x_t(r+a(t),-r+a(t))}$ (with respect to $\Omega_t$).

The first modification we need to introduce is composing all curves  used in the proof of theorem 
 \ref{thm:equiv} and defined in a neighborhood of $0\in \mathbb{C}$,  with the  diffeomorphism
\[(x,y)\mapsto (x+a(y),y).\]

The second difference is that from the very beginning our parallel transport here
 is with respect to a family of symplectic forms $\Omega_t$, and with $\Omega_0$
 of type $(1,1)$ at the origin. This situation is not quite new since in the  proof of theorem  \ref{thm:equiv}
we already needed to interpolate symplectic forms (although at a later stage).

Hence we conclude that for $r$ small enough the fiber $W_{-r}$ is  equivalent to Lagrangian surgery
 (and hence by theorem \ref{thm:equiv} generalized Dehn surgery)
along a framed Lagrangian sphere $L$; the Lagrangian sphere  are the points in $W_r$ which parallel transport over $x_0(r,0)$
 sends to the critical point $c$.
\end{proof}

\begin{remark} Theorem \ref{thm:handle} is rather natural in view of the results for contact manifolds
 in \cite{Pr02}.
\end{remark}

\subsection{Further directions.}

In this paper we have shown that 2-calibrated foliations are a wide enough class 
of codimension one foliations, and not surprisingly,
techniques from symplectic geometry are well suited to their study. 
 We would like to finish by discussing a couple of questions that we were not able to answer.

 Theorem \ref{thm:leafhomeo} shows that our embedded 3-dimensional taut foliations capture the
 leaf space of $\mathcal{F}$. What it would be interesting to know is whether they capture the full 
transverse geometry, i.e. the holonomy groupoid.

A remarkable property of 3-dimensional taut foliations is that transverse loops are never nullhomotopic.
 The proof of this fact uses that the universal cover of the 3-manifold is $\mathbb{R}^3$, 
a property which does not extend to manifolds supporting a 2-calibrated foliation.
 We know no examples 2-calibrated foliations on simply connected manifolds:  in \cite{IM03} it was shown that 
 normal connected sum could be used to construct 5-dimensional simply connected regular Poisson
manifolds with codimension one leaves. Those methods,  however,
cannot be used to construct simply connected 2-calibrated for conditions
 in theorem \ref{thm:calibnorsum} are not fulfilled. It has been recently shown that
Lawson's foliation on $S^5$ is the symplectic foliation of a Poisson structure \cite{Mi11}. However,
this Poisson structure does not admit a 2-calibration because Lawson's foliation is not taut (the compact leaf 
would make any transverse loop non-trivial in homology). 

We conjecture that any transverse loop in a 2-calibrated foliation is not nullhomotopic.

\section{Appendix: Legendrian surgery, open book decompositions and generalized Dehn
surgery}\label{ssec:Legendrian}

Let $(M,\xi)$ be an exact contact manifold and let $\alpha$ be a contact 1-form 
defining $\xi$  ($\xi=\mathrm{Ker}\alpha$). Recall that an open book decomposition for $M$ is given by a pair
$(K,\theta)$ such that
\begin{itemize}
 \item   $K$ is a codimension 2 submanifold with trivial normal bundle, referred to as
the binding;
 \item  $\theta\colon M\backslash K\rightarrow S^1$ is  a
fibration that  in a trivialization  $D^2\times K$ of a neighborhood of $K$ is the angular coordinate.
 \end{itemize}

 Let  $F$ denote the closure of any fiber of  $\theta$. The first return map  associated
  to a suitable  lift of  $\tfrac{\partial}{\partial \theta}$ to $M\backslash K$,
   defines a diffeomorphism of $F$  supported away
    from a neighborhood of the boundary $\partial{F}=K$. Up to diffeomorphism $M$ can be recovered out of $F$ and the
first return map.

The following discussion  is mostly taken from \cite{GM02}:

\begin{definition} The contact structure $\xi$ is supported by an open
 book decomposition $(K,\theta)$ if for a choice of contact form $\alpha$ defining $\xi$ we
have:
\begin{itemize}
\item $\alpha$ restricts to $K$ to a contact form.
\item $d\alpha$ restricts to each fiber of $\theta$ to an exact symplectic structure.
\item The orientation of $K$ as the boundary of each symplectic leaf matches the
natural orientation induced by the contact form.
\end{itemize}

The form $\alpha$ is said to be adapted to the open book decomposition
$(K,\theta)$.
\end{definition}

In what follows we are going to discuss contact structures and cosymplectic foliations on a given manifold.
Since we have been using the notion of Reeb vector field for cosymplectic foliations,  we
refer to contact Reeb vector fields when discussing contact structures.

Given a contact form  $\alpha$  adapted to $(K,\theta)$, it is possible to scale it away from $K$ to a contact 1-form $\alpha'$
such that the flow along its contact Reeb vector field defines a compactly supported  first
return map $\varphi\in \mathrm{Symp}(\mathrm{int}F,d\alpha')$ \cite{GM02}.

The isotopy class of $(M,\xi)$ is totally determined by any
open book decomposition supporting it \cite{Gi02,GM02}. More precisely, the relevant structure in the open book decomposition
is  the completion of the structure of exact symplectic manifold convex at infinity  of
the exact symplectic fiber  $(\mathrm{int}F,d\alpha)$ (or $(\mathrm{int}F,d\alpha')$), together
with the first return symplectomorphism supported inside $\mathrm{int}F$.

The previous characterization becomes very important in light of the following theorem:
\begin{theorem}\cite{Gi02,GM02}\label{thm:contbook} For every exact contact manifold
 $(M,\xi)$ and
any contact form defining $\alpha$, there exist an open book decomposition  $(K,\theta)$
supporting $\xi$  such that $\alpha$ is adapted to it.
\end{theorem}

Let $\alpha$ be a contact form on $M$ 
 adapted to the open book decomposition  $(K,\theta)$ and let  $L$
be a parametrized Legendrian sphere which is contained in a fiber of $\theta$, and hence
it becomes Lagrangian for the symplectic structure $d\alpha$ on the fiber.
 
Observe that  away from the binding $K$, the open book decomposition defines a
2-calibrated foliation $(M\backslash K,\mathcal{F}_\theta,d\alpha)$, with $\mathcal{F}_\theta=\mathrm{Ker}d\theta$, 
which is a symplectic mapping torus associated to the symplectomorphism $\varphi$ supported in $\mathrm{int}F$. 
 Generalized Dehn surgery along $L$  produces a new symplectic mapping torus
 with return map $ \varphi\circ\tau$, where $\tau$ is a generalized
Dehn twist along $L$. Because the symplectic leaf is the same and  the return map is still compactly
 supported,  the symplectic mapping torus is in fact the open book decomposition of a unique contact manifold 
(up to isotopy). In \cite{GM02} it has been announced that this contact manifold is $(M^L,\alpha^L)$ the result 
of performing Legendrian surgery along $L$ \cite{We91}.

 The ideas developed relating Lagrangian surgery and generalized Dehn
 surgery allow us to give a very natural proof of this result. The key step is the following theorem.

\begin{theorem}\label{thm:leg} Let $L\subset (M,\alpha)$ be a parametrized Legendrian sphere 
in a contact manifold and let $(M^L,\alpha^L)$ be the contact
manifold obtained by Legendrian surgery along $L$. Suppose that $\alpha$ is adapted to the open
 book $(K,\theta)$ and that $L$ is contained in a fiber of $\theta$.
Then  given $V$ any small enough neighborhood of $L$ with empty intersection
with the binding $K$, there exists an isotopy $\Psi_s\colon M\rightarrow M$,
$s\in [0,1]$, starting at the identity with the following properties:
\begin{itemize}
\item $\Psi_s$ is supported inside $V$ and tangent to the identity at $L$.
\item  $(M\backslash K,\mathcal{F}_{\theta_s},d\alpha)$, with
$\mathcal{F}_{\theta_s}:={\Psi_s}_{*}\mathcal{F}_\theta$,
 is a  2-calibrated foliation  and thus an open book decomposition  $(K,{\Psi_s}_*\theta)$ of $M$ 
to which the contact form $\alpha$ is adapted.
\item Let $(M^L\backslash
K,\mathcal{F}_{\theta_1^{L}},d\alpha^{L})$ be the result of
performing generalized Dehn surgery on $(M\backslash
K,\mathcal{F}_{\theta_1},d\alpha)$ along the parametrized Lagrangian sphere $L$. Then  $(M^L\backslash
K,\mathcal{F}_{\theta_1^{L}},d\alpha^{L})$ is an open book decomposition  $(K,\theta_1^{L})$ for  $M^L$ and
 the contact form $\alpha^{L}\in \Omega^1(M^{L})$ is adapted to $(K,\theta_1^{L})$.
\end{itemize}
\end{theorem}

\begin{proof}
We are going to recall Weinstein's definition of Legendrian surgery using a symplectic
cobordisms and a Liouville vector field transverse to the boundary. Actually, we will modify  the original choices to make them
 compatible with our setup for Lagrangian surgery, or by theorem
\ref{thm:equiv} with the setup for generalized Dehn surgery.

Recall that a boundary component of a symplectic manifold $(Z,\Omega)$ (of dimension bigger than 2)
 endowed with a Liouville vector field $Y$ 
is said to be  convex (respectively concave) if $Y$ is   outward (respectively inward) pointing. 

We consider $(M\times [-1,1],d(e^v\alpha))$, which is a subset of the symplectization of
$(M,\alpha)$. The tuple $(M\times [-1,1],d(e^v\alpha), \tfrac{\partial}{\partial v}, M\times
\{0\}, L\times\{0\})$ is an isotropic setup in  the language of
Weinstein (\cite{We91}, section Neighborhoods of isotropic submanifolds). Note that $\{1\}\times M$
(respectively $\{-1\}\times M$) is a convex (respectively) concave boundary component 
(beware that the  notion of  Liouville vector field we
 use is opposite to Weinstein's, for we require the flow of the vector field to expand  the symplectic form exponentially).

The second  isotropic setup is the one of the $(n+1)$-handle to be
attached, which is the one described in  \cite{We91}, section Standard Handle, up to the following change. Unlike Weinstein does,
 we are going to glue the convex end of
 $(M\times [-1,1],d(e^v\alpha), \tfrac{\partial}{\partial v}, M\times
\{0\}, L\times\{0\})$ to the concave end of the symplectic $(n+1)$-handle; the reason is that in our definition
of Lagrangian surgery, we glued  the symplectic $(n+1)$-handle along the hypersurface $H_{2,r}$ 
 where the symplectic vector field points inward. 
 For this reason we also define a different Liouville vector field in the $(n+1)$-handle.
 We  use the notation introduced in
\ref{ssec:lagsurg}.

 The symplectic form is the standard one  $\Omega_{\mathbb{R}^{2n+2}}$. We consider the function
 \[q=\sum_{i=1}^{n+1}x_i^2-2y_i^2,\]
whose negative gradient with respect to the Euclidean metric
 \[E=-2x^1\frac{\partial}{\partial x^1}+4y^1\frac{\partial}{\partial y^1}-\dots
-2x^{n+1}\frac{\partial}{\partial
x^{n+1}}+4y^{n+1}\frac{\partial}{\partial y^{n+1}}\]
 is  a Liouville vector field.

 For each $r>0$ we consider the fiber $q_r$, which  contains the
Lagrangian sphere  $\Sigma_r$  described in lemma \ref{lem:stb} using $Y_2$  the Hamiltonian vector field of $-\mathrm{Re}h$
with respect to  $\Omega_{\mathbb{R}^{2n+2}}$. Notice that $dq(Y_2)<0$ and therefore
$Y_2$ is transverse to the level hypersurfaces  $q_r$. Since $Y_2$ and $E$ coincide at $\Sigma_r$, it follows that
the sphere $\Sigma_r$ is also Legendrian with respect to the contact form
$\alpha_E:=i_E\Omega_{\mathbb{R}^{2n+2}}$ on $q_r$. Moreover, at points of $\Sigma_r\subset q_r$
 the contact distribution and the cosymplectic distribution coincide.

 Let $V_r(\epsilon)$ be a
tubular neighborhood or radius $\epsilon>0$ of $\Sigma_r$ inside
$q_r$ with respect to the Euclidean metric.
We claim that for any $\epsilon'>0$, $\epsilon >\epsilon'$, we have
$f_r\in C^{\infty}(V_r(\epsilon)\backslash \Sigma_r,\mathbb{R}^+)$ a cut-off
function with compact support and with the following two properties:
\begin{itemize}
\item $\varPhi_1^{f_rY_2}(V_r(\epsilon')\backslash \Sigma_{r})\subset
q_{-2r}$ (note that $q_{-2r}$ contains the Lagrangian sphere $\Sigma_{-r}$).
\item $\varPhi_1^{f_rY_2}(V_r(\epsilon))$ is transverse to $E$.
\end{itemize}
Assuming the claim, we define the hypersurface
\[H^{L}_r:=\varPhi_1^{f_rY_2}(V_r(\epsilon)\backslash
\Sigma_{r})\cup \Sigma_{-r}.\]
By assumption  the Liouville vector field $E$ is
transverse to $H^{L}_r$, and thus the hypersurface inherits an
exact contact structure $\alpha_{E}$ by restricting $i_E\Omega_{\mathbb{R}^{2n+2}}$.

The second isotropic setup is the following: the symplectic $(n+1)$-handle is the compact region bounded by
$H^{L}_r$ and $V_r(\epsilon)$ endowed with the standard symplectic form; the Liouville vector field is $E$;
the hypersurface is $V_r(\epsilon)$, which is concave; the parametrized Legendrian sphere is $\Sigma_r$.

The symplectic morphism $\psi$ that gives rise to the symplectic elementary cobordism 
  (\cite{We91}, proposition 4.2, whose replacement for Lagrangian
 surgery is lemma \ref{lem:cosympglue}), sends $(V_r(\epsilon),\Sigma_r,\alpha_E)$ to
$(\nu(L),L,\alpha)$,  and therefore we can consider
$(V_r(\epsilon),\Sigma_r,\alpha_E)$ as a subset of
$(M,\alpha)$. Then
 \[M^{L}:=H^{L}_r\cup (M\backslash 
V_r(\epsilon))\] carries and obvious contact form $\alpha^{L}$
which extends $(M\backslash V_r(\epsilon),\alpha)$.

The data for Legendrian surgery has been chosen to be compatible with Lagrangian surgery:
 both  $H^{L}_r$ and $V_r(\epsilon)$ are
transverse to $Y_2$  and therefore they inherit 2-calibrated foliations
$(H^{L}_r,\mathcal{F}_r^L,\omega^L_r)$ and $(V_r(\epsilon),\mathcal{F}_r,d\alpha)$.
Theorem \ref{thm:equiv} easily implies that
$(H^{L}_r,\mathcal{F}_r^L,\omega^L_r)$ is the result of generalized Dehn surgery along  
$\Sigma_r\subset (V_r(\epsilon),\mathcal{F}_r,d\alpha)$.

On $V_r(\epsilon)$ we have two structures of 2-calibrated foliation,
$(\mathcal{F}_r,d\alpha)$ and $(\mathcal{F},d\alpha)$. The reason is that $\psi$ preserves
contact forms and hence contact Reeb vector fields, but it does not preserve the 1-forms
defining the cosymplectic foliations (or their associated Reeb vector fields). However, at $\Sigma_r$ the Liouville
and Hamiltonian vector field coincide, and this implies that at points in $L$ the contact distribution is tangent
to $\mathcal{F}_r$. In particular the contact Reeb vector field for $\alpha$ is transverse to 
$\mathcal{F}_r$ near $L$. It is also transverse to $\mathcal{F}$ because $\alpha$ is adapted to the open book.
Therefore we can use
the trajectories of the contact Reeb vector field, to construct an isotopy $\Psi_s$ tangent to the identity at $L$ and supported 
inside $V$  a small neighborhood of $\Sigma_r$ contained in $V_r(\epsilon)$.

The claim about the existence of the function $f_r$ is easily proved
when $n=1$ by inspecting the trajectories of  $E$ and $Y_2$. The
general case can be reduced to the previous one: each point
$(x_1,y_1,\dots,x_{n+1},y_{n+1})$ in $\mathbb{C}^{n+1}$ and away from the
union of stable and unstable manifolds (these
are the same for  both Morse functions $\mathrm{Re}h$ and $q$),
determines $[x_1:\cdots :x_{n+1}],[y_1:\cdots :y_{n+1}]$ a point in
$\mathbb{R}\mathbb{P}^{n}\times \mathbb{R}\mathbb{P}^{n}$, which
gives rise to two lines in $\mathbb{R}^{n+1}$ and
$i\mathbb{R}^{n+1}$ respectively. These lines span a plane in
$\mathbb{C}^{n+1}=\mathbb{R}\oplus i\mathbb{R}^{n+1}$. Each plane in the family is
preserved by the flow of $E$ and $Y_2$; moreover, the flows restrict to the planes to
the flows of the 1-dimensional case. From this observation the claim
follows easily.
\end{proof}

Theorem \ref{thm:leg} provides an isotopy $\Psi_s$ supported away from $K$, so that $\alpha$ is adapted to the 1-parameter 
family of open book decompositions $(K,\Psi_s\theta)$. Therefore  we can identify
the symplectic fiber $F$ and symplectic monodromy $\varphi\in
\mathrm{Symp}(\mathrm{int}F,d\alpha)$ of $(K,\theta)$  with those of $(K,\Psi_1\theta)$ (again following the contact Reeb flow).
Hence point 3 in theorem \ref{thm:leg} asserts that  
$(M^{L},\alpha^{L})$ is adapted to an open book decomposition with 
  the same symplectic leaf $(F,d\alpha)$ and monodromy $\varphi\circ \tau \in
\mathrm{Symp}(\mathrm{int}F,d\alpha)$, which is exactly what we wanted to prove.

\begin{remark}
If we attach the convex end of the symplectic handle to the concave end of the
 symplectization, we get the contact manifold $(M^{L^-},\alpha^{L^-})$. It is easy to see that 
 $\alpha^{L^-}$ is  adapted to an open book decomposition whose monodromy is $\varphi\circ \tau^{-1}$.

Observe that proposition 6.1 in \cite{DK75}  implies that in  dimensions 5 and 13
 the manifolds $M^L$ and $M^{L^-}$ are diffeomorphic. In \cite{KN}, section 3, it is shown
 that there are instances (coming from Brieskorn manifolds) in which
 $(M^L,\alpha^L)$ and $(M^{L^-},\alpha^{L^-})$ are not contactomorphic,
  and hence the authors can deduce that $\tau^2$ is not isotopic to
  the identity in
$\mathrm{Symp}^{\mathrm{comp}}(T^*S^6,d\alpha_{\mathrm{can}})$,
   a result  already proved by Seidel for $n=2$ \cite{Se98}; similar results
    are also drawn for powers of the Dehn twists known to be isotopic to
    the identity in  $\mathrm{Diff}^{\mathrm{comp}}(T(\lambda))$, for all $n$
even.
\end{remark}

\begin{remark}
For any contact form $\alpha$ on $M$ representing the given contact
structure $\xi$ and $L$ a Legendrian submanifold, Giroux and Mohsen announce
\cite{GM02} the existence
of relative open book decompositions, meaning that $\alpha$ is adapted to the open book decomposition and $L$ is contained
 in a fiber.

The interested reader familiar with approximately holomorphic
geometry  \cite{Do96} and its version for contact manifolds \cite{IMP00,Pr02},
  can write a proof along the following lines: the open book decomposition is the
  result of pulling back the canonical open book decomposition of $\mathbb{C}$ by an
   approximately holomorphic function. To make sure the binding does not contain
$L$,
    we use reference sections supported near $L$ which achieve the value 1 when
restricted
     to $L$; they come from an explicit formula once we identify a tubular
neighborhood
      of $L$ with a tubular neighborhood of the zero section of the first jet bundle
with its canonical contact structure $(\mathcal{J}^1L,\alpha_{\mathrm{can}})$.  It is necessary to further add
perturbations whose restriction to $L$ attain real values: they are such that its
   restriction to $T^*L\times \{0\}\subset \mathcal{J}^1L$ are small
real multiples of reference sections equivariant with respect to 
         the involution on $(\mathcal{J}^1,\alpha_{\mathrm{can}})$
         which reverses the sign of the fiber and conjugation on $\mathbb{C}$ (this construction is analogous to the
content of the remark after lemma 3 in \cite{AM}).

Therefore we conclude that Lagrangian surgery includes Legendrian
surgery, for we can bypass the latter by choosing appropriate
compatible open book decompositions and then performing Lagrangian
surgery. According to theorem \ref{thm:equiv} we can even claim that
generalized Dehn surgery contains Legendrian surgery, and  forget
about the cobordisms.

Actually, the reason why generalized Dehn surgeries for different open book decompositions supporting
the contact structure give the same contact manifold, is because there is a
 contact surgery behind. Now consider $(L,\chi)$ where $L$ is a Legendrian
  submanifold of $(M,\alpha)$ and $\chi\in
\mathrm{Symp}^{\mathrm{comp}}(T^*L,d\alpha_{\mathrm{can}})$.
Let us   take any open book decomposition relative to $L$ and such that $\alpha$ is adapted to it,
 and  consider the new manifold $M^L$ associated to the open book decomposition with symplectic
   monodromy $ \varphi\circ \chi$. It is clear that the diffeomorphism type of
the manifold
    does not depend on the open book decomposition, but it is not clear  whether in general
the contact
     structure depends on the choice of open book decomposition. In either case, it would be
an interesting situation because it would give either a new contact surgery -possibly  a
Legendrian surgery
       based on a block different from a symplectic handle- or different
contact structures.
\end{remark}

%
%
%

\end{document}